\documentstyle[11pt]{article}
\input psfig
\def\IMSmarkvadjust{0 pt}
\def\IMSmarkhadjust{0 pt}
\def\IMSmarkhpadding{0 pt}
\def\IMSpubltext{Published in modified form:}
\def\SBIMSMark#1#2#3{
 \font\SBF=cmss10 at 10 true pt
 \font\SBI=cmssi10 at 10 true pt
 \setbox0=\hbox{\SBF \hbox to \IMSmarkhpadding{\relax}
                Stony Brook IMS Preprint \##1}
 \setbox2=\hbox to \wd0{\hfil \SBI #2}
 \setbox4=\hbox to \wd0{\hfil \SBI #3}
 \setbox6=\hbox to \wd0{\hss
             \vbox{\hsize=\wd0 \parskip=0pt \baselineskip=10 true pt
                   \copy0 \break%
                   \copy2 \break%
                   \copy4 \break}}
 \dimen0=\ht6   \advance\dimen0 by \vsize \advance\dimen0 by 8 true pt
                \advance\dimen0 by -\pagetotal
	        \advance\dimen0 by \IMSmarkvadjust
 \dimen2=\hsize \advance\dimen2 by .25 true in
	        \advance\dimen2 by \IMSmarkhadjust

  \openin2=publishd.tex
  \ifeof2\setbox0=\hbox to 0pt{}
  \else 
     \setbox0=\hbox to 3.1 true in{
                \vbox to \ht6{\hsize=3 true in \parskip=0pt  \noindent  
                {\SBI \IMSpubltext}\hfil\break
                {\it  J. Funct. Anal.}~{\bf 143} (1997), 309--336
 
                \vfill}}
  \fi
  \closein2
  \ht0=0pt \dp0=0pt
 \ht6=0pt \dp6=0pt
 \setbox8=\vbox to \dimen0{\vfill \hbox to \dimen2{\copy0 \hss \copy6}}
 \ht8=0pt \dp8=0pt \wd8=0pt
 \copy8
 \message{*** Stony Brook IMS Preprint #1, #2. #3 ***}
}

\renewcommand{\baselinestretch}{1.3}

\newtheorem {th}{Theorem}[section]
\newtheorem {lem}[th]{Lemma}

\newtheorem {pr}[th]{Proposition}
\newtheorem {cor}[th]{Corollary}

\newcommand{\ben}{\begin{enumerate}}
\newcommand{\een}{\end{enumerate}}

\newcommand{\bea}{\begin{eqnarray}}
\newcommand{\ba}{\begin{array}}
\newcommand{\bean}{\begin{eqnarray*}}
\newcommand{\ea}{\end{array}}
\newcommand{\eea}{\end{eqnarray}}
\newcommand{\eean}{\end{eqnarray*}}
\newcommand{\be}{\begin{equation}}
\newcommand{\ee}{\end{equation}}

\newcommand{\Q}{{\rm I\!\!\!Q}}
\alph{enumii}
\roman{enumiii}

\def\Cox{\fbox{\hphantom{} \vphantom{x}}}

\def\sf{\sigma\mbox{-field}}

\def\eps{\epsilon}
\def\E{{\bf{E}}}
\def\P{{\bf{P}}}
\def\N{\hbox{I\kern-.2em\hbox{N}}}
\def\R{\hbox{I\kern-.2em\hbox{R}}}

\def\C{{\cal{C}}}

\newcommand{\ZZ}{Z\!\!\!Z}
\def\D{{\cal{D}}}

\def\F{{\cal{F}}}

\def\cA{{\cal{A}}}

\def\|{\, | \, }
\def\cond{\, \Big| \, }
\def\one{\mbox{\bf  1}}
\def\0{\hat{0}}
\def\1{\hat{1}}

\def\xx{{\bf x}}

\def\trunc{\mbox{ \rm ~trunc }}
       
\def\cen{\mbox{\rm center}}
\def\diam{\mbox{\rm diam}}
\def\Wall{\mbox{\rm Wall}}
\def\pare{\mbox{\rm parent}}
\def\dist{\mbox{\rm dist}}

\def\gs{||G_*||}
\def\clol{c_1}
\def\cllo{c_2}
\def\clo{c_3}
\def\cl{c_4}
\def\cll{c_5}
\def\clll{c_0}
\def\cbeta{c_6}
\def\cloop{c_7}
\def\cfast{c_8}
\def\clool{c_9}
\def\clolo{c_{10}}
\def\fr{\mbox{\rm \small frontier}}
\def\frK{\mbox {\rm  \small frontier}(K)}
\def\frKB{\mbox{\rm  \small frontier}(B[0,1])}
\def\frKBt{\mbox{\rm  \small frontier}(B[0,t])}
\def\frKBtau{\mbox{\rm  \small frontier}(B[0,\tau])}
\def\esssup{\mbox {\rm ess sup}\;}
\def\trunc{\mbox{\rm trunc}}

\def\bcorG{\mbox{{\bf core}$(G,\eta)$ }} 
\def\corG{\mbox{{\rm  core}$(G,\eta)$ } }
\def\core{\mbox{\rm core}}
\def\crG0{\mbox{{\rm  core}$(G,\eta_0) \, $}}  
\def\wtg{\widetilde \gamma}
\def\rad{\mbox{\rm rad}}
\def\tauc{\tau_{_{\rm C}}}
\def\taucG{\tau_{_{{\rm C}, \partial G}}}
\def\taucJ{\tau_{_{{\rm C}, J}}}
\def\tauG{\tau_{_{ \partial G}}}
\def\texp{\widetilde \tau_{\exp}}
\def\complex{\mathop{\raise .45ex\hbox{${\bf\scriptstyle{|}}$}
    \kern -0.40em {\rm \textstyle{C}}}\nolimits}

\def\Ubar{\overline{U}}
\def\Bbridge{B_{\rm br}}
\begin{document}

\begin{titlepage}
\begin{center}
{\large \bf The Dimension of the Brownian Frontier is Greater Than 1.}
\end{center}
\begin{center}
{\sc Christopher J. Bishop}\footnote{Supported in part by
NSF grant \# DMS 9204092 and by an Alfred P.
Sloan Foundation  Fellowship.},
%
\, {\sc Peter W. Jones}\footnote{Research partially 
 supported by NSF grant \# DMS 9213595.}, \\
{\sc Robin Pemantle}\footnote{Research supported in part by
National Science Foundation grant \# DMS 9300191, by a Sloan Foundation
   Fellowship, and by a Presidential Faculty
Fellowship.}
%
and {\sc Yuval Peres}\footnote{Research partially supported by NSF grant
\# DMS-9404391.}
\end{center}
\begin{center}
\it State University of New York, Yale University, \\
University of 
Wisconsin, and University of California.
\end{center} \rm
\vspace{.3in}
\centerline{\bf Abstract}
Consider a planar Brownian motion run for finite time. 
 The {\em frontier} or ``outer boundary'' of the path is the boundary 
of the unbounded component of
the complement.  Burdzy (1989) showed that the frontier has
infinite length.  We improve this by showing that the Hausdorff
dimension of the frontier is strictly greater than 1.
(It has been conjectured that the Brownian frontier has dimension
 $4/3$, but this is still open.)  The
proof  uses  Jones's Traveling Salesman Theorem
and  a self-similar tiling of the plane 
by fractal tiles known as Gosper Islands.  

\vspace{1in}

\SBIMSMark{1995/9}{June 1995}{}
\end{titlepage}

\section{Introduction} \label{sec 1}

Let $K$ be any compact, connected set in the plane.
The complement of $K$ has one unbounded component and its 
topological boundary is called the {\em frontier} of $K$,
denoted $\frK$.  The example we are most interested in is when $K$ 
is the range of a planar Brownian motion run for a finite time
(see Figure~\ref{outer}).
In this case, Mandelbrot (1982) conjectured that the Hausdorff
dimension $\dim (\frK)$ is $4/3$.  Rigorously, the best proven upper
bound on the dimension is $3/2 - 1/(4\pi^2) \approx 1.475$ by
Burdzy and Lawler (1990).  Burdzy (1989) proved that
$\frK$ has infinite length; our main result improves this
to a strict dimension inequality:
\begin{figure}
\vbox{ \centerline{ \hbox{
\psfig{figure=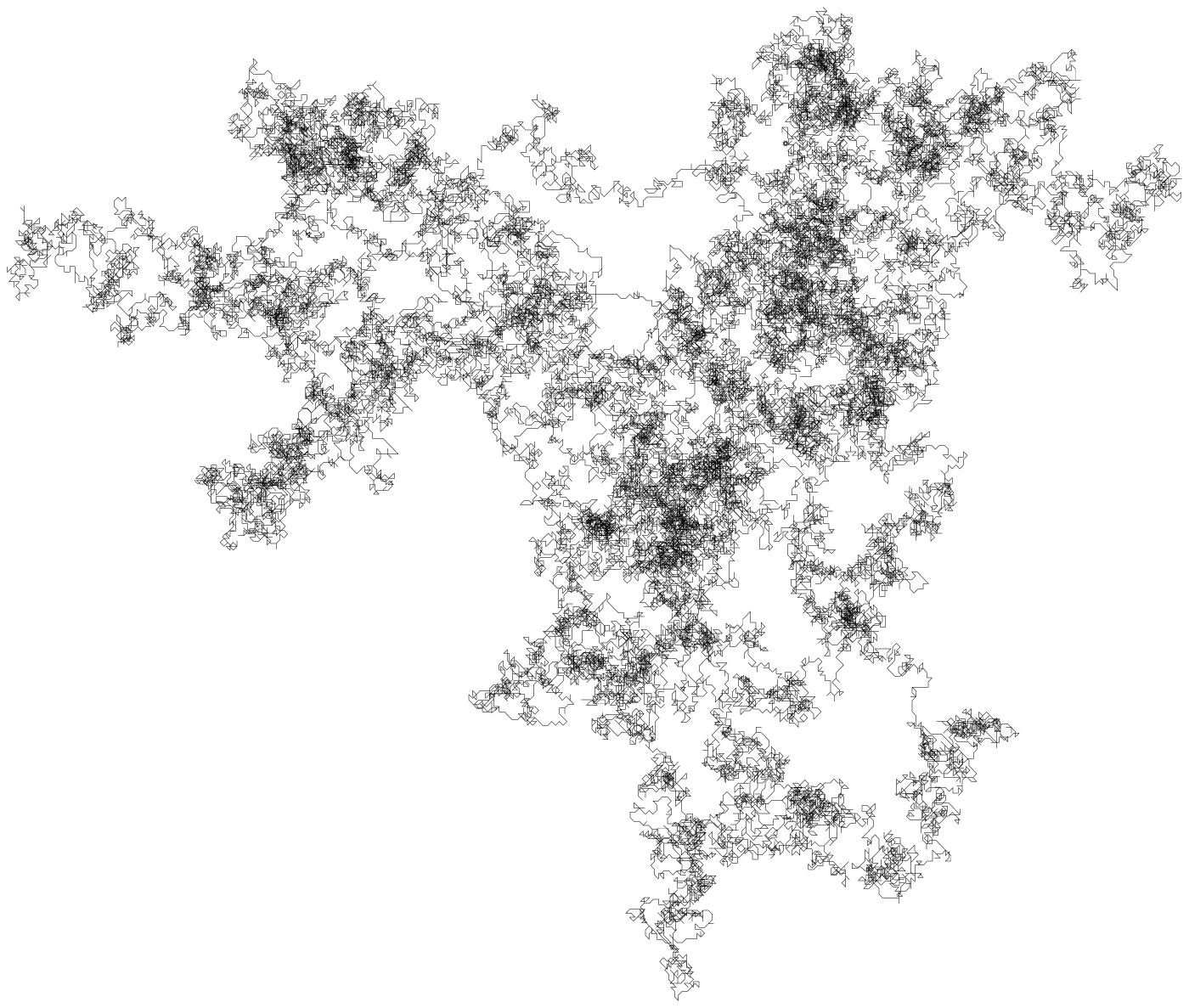,height=2.5in}
$\hphantom{xxx}$
\psfig{figure=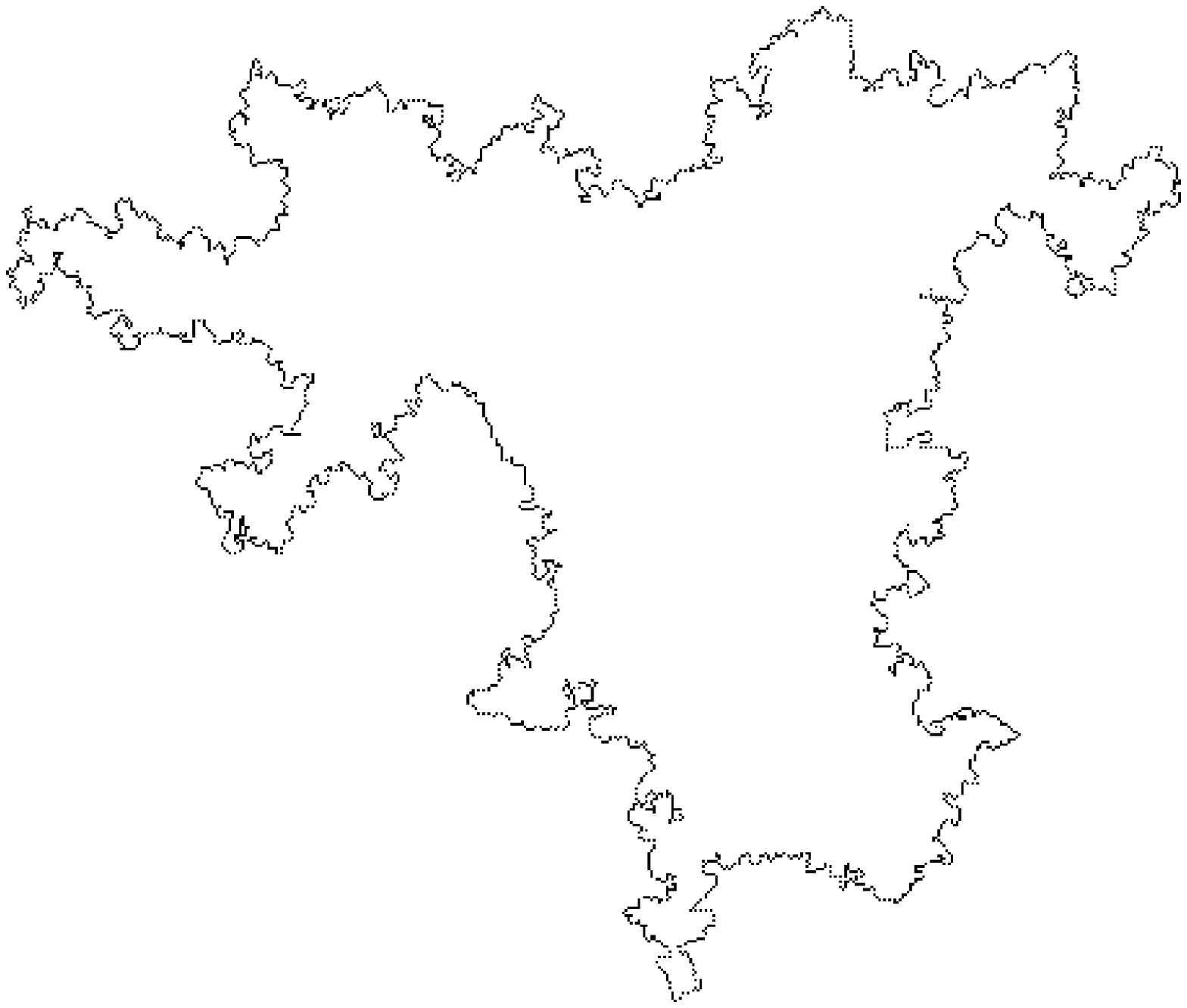,height=2.5in}
}}}
 \caption{ \label{outer} A Brownian path and its frontier}
\end{figure}

\begin{th} \label{th main}
Let $B[0,t]$ denote the range of a planar Brownian motion, 
 run until time $t>0$.  There is an $\eps > 0$ such that with probability 1,
The Hausdorff dimension $ \, \dim (\frKB)$ is at least $1 + \eps$. 
Moreover, with probability 1,
$$
\inf_{t>0} \, \inf_V \, \dim \Big(\frKBt \cap V \Big) \, \geq \, 1 + \eps \, ,
$$   
where the inner infimum is over all open sets $V$ that intersect $\, \frKBt$.
\end{th}

\noindent{\em Remarks:}   The uniformity in $t$ implies that 
$ \, \dim (\frKBtau) \geq 1 + \eps$ almost surely for any
positive random variable $\tau$ (which may depend on the Brownian motion).
We also note that
 our proof shows that the frontier can be replaced in the statement of
the theorem by the boundary of any connected component of 
the complement $B[0,t]^c$.
  (One can also infer this from the statement of the theorem
 by using  conformal invariance of Brownian motion).
As explained at the end of Section \ref{sec 6},
The result also extends to the frontier of the planar Brownian
bridge (which is a closed Jordan curve by Burdzy and Lawler (1990)).

Bishop and Jones (1994) proved that if a compact, connected set
is ``uniformly wiggly at all scales'', then it has dimension
strictly greater than 1.  Here we adapt this to a stochastic setting
in which the set is likely to be wiggly at each scale, given
the behavior at previous scales.  The difficulty is in handling
statistical dependence.

\noindent{\bf Definitions:} Let $G$ be a compact set in the plane
with complement $ G^c$, and let
 $\eta>0$. Denote by \bcorG the set 
$\{ z \in G  : \dist (z,  G^c) > \eta \cdot \diam(G) \}$.
Say that the compact set $K \, $ \mbox{\bf $\, \eta$-surrounds} $G$ if $K$
topologically separates \corG from $ G^c$,
i.e., if  \corG is disjoint from the unbounded component of $(K \cap G)^c$.

\begin{figure}
\centerline{ \psfig{figure=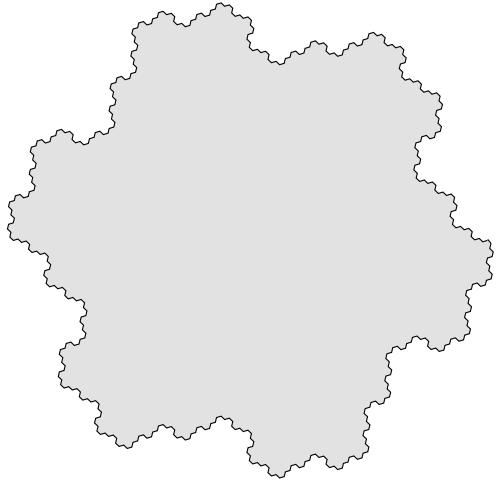,height=2.5in} }
 \caption{ \label{gen4} The Gosper island }
\end{figure}
\begin{th} \label{th any K}
Let $G_0$ be the {\bf Gosper Island}, defined in the next section
 and illustrated in
Figure~\ref{gen4}. There exists an  absolute
constant $\eta_0 > 0$ with the following property.
 Suppose that $ \clll >0$, and $K$ is a random compact connected
subset of the plane such that for 
 all homothetic
images $G=z+r G_0$ of $G_0$ with $r \in (0, 1)$ and $z$ in the plane:  
\begin{equation} \label{eq gosper-wiggly}
   \P \Big [ K \, \mbox{ \rm  $\eta_0$-surrounds } G 
             \, \mbox { \rm \underbar{or} } \; \, K \cap \crG0 = \emptyset
      \, \cond \, \sigma(K \setminus G^{\circ} ) \Big ] > \clll,
\end{equation}
where the conditioning is on the $\sigma$-field generated by
the random set $K$ outside the interior $G^{\circ}$ of $ G$. 
Then there is an $\, \eps > 0 $,  depending only on $\clll$, such that
$$
\dim (\frK) \geq 1 + \eps
$$
with probability 1.  More generally, $\dim (\frK \cap V) \geq 1 + \eps$
for any open $V$ intersecting $\frK$.  
\end{th}
\begin{figure}
\centerline{ \psfig{figure=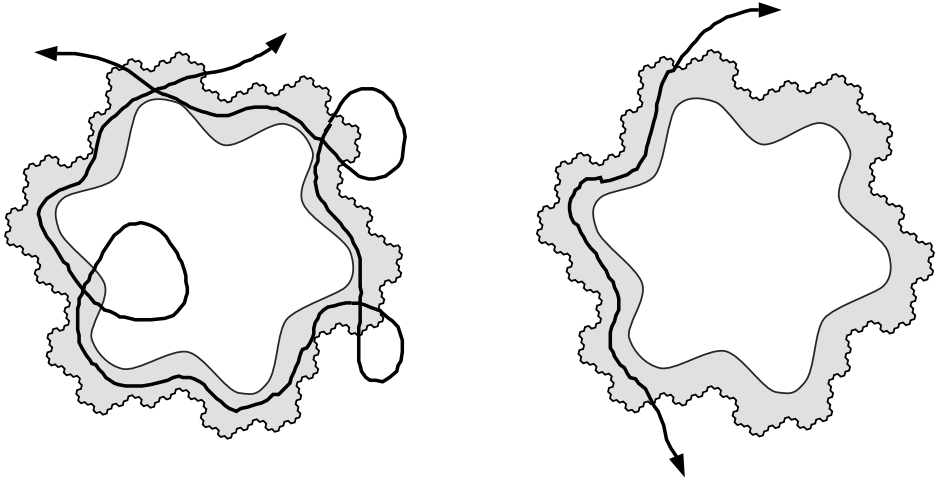,height=2.5in} }
 \caption{ \label{Surr+Miss} A Brownian motion which surrounds the core
and one which misses it. }
\end{figure}

\noindent{\sc Remarks: } { \bf 1.} In fact, the proof in Section $ 3 \,$
 shows that
with probability 1, for any connected component  $\Omega$
of $K^c$ and any open $V$ intersecting $\partial \Omega$,
there is a John domain $\Omega_{\rm J} \subset \Omega$ with 
closure $\overline{\Omega_{\rm J}}$  contained in
$V$, such that 
$\dim (\partial \Omega \cap \partial \Omega_{\rm J}) \geq  1 + \eps$.    
 \newline {\bf 2.} The constant $\eta_0$ 
will be chosen in the next section to ensure that no ``macroscopic''
line segment can be wholly contained within a $2 \eta_0$-neighborhood
of the Gosper Island's boundary $\partial G_0$.

The appearance of the Gosper Island might seem strange at this point, 
but is explained as follows.  The hypothesis on $K$ that
guarantees ``wiggliness'' should be local to  handle dependence (thus
it must hold inside each $G$ conditioned on $K \cap G^c$).
If $ K \, $  $\, \eta_0$-surrounds $G$, 
     \underbar{or} $\;  K \cap \crG0 = \emptyset$,
then $\frK$ cannot intersect \crG0.   
Having thus controlled $\frK$ inside $G$, away from the
boundary of $G$, we must worry about how $\frK$ behaves
near boundaries of cells $G$, as these run over a 
partition of the plane.  If a small neighborhood of the union 
of the boundaries of  cells $G$ of a fixed size
contains no  straight line segments of length comparable to
$\, \diam \, (G) \,$, then
no significant flatness can be introduced near cell boundaries.
To apply the argument with the same constants on every scale,
we need a self-similar tiling where tile boundaries have
no straight portions; the Gosper Island yields such a tiling.  

Proving  Theorem \ref{th any K} is the main effort of the paper and
is organized as follows.  Section~\ref{sec 2} summarizes
notation and useful facts about the Gosper Island.  
We also discuss the notion of a Whitney decomposition with respect 
to these tiles.  Section~\ref{sec 3} constructs a random tree 
of Whitney tiles for $K$ and reduces Theorem \ref{th any K} to a
lower bound on the expected growth rate of the tree, via some general
propositions on random trees.
In Section~\ref{sec 4} we state a variant of 
Jones's Traveling Salesman Theorem adapted to 
the current setting.  In Section~\ref{sec 5} this theorem is used
to derive the required lower bound on the expected growth rate
of the ``Whitney tree'' mentioned above, which 
 then finishes the proof of Theorem~\ref{th any K}.  In
Section~\ref{sec 6} we verify that  the range of planar Brownian motion,
killed at an independent exponential time,
satisfies the hypothesis of Theorem~\ref{th any K};
this easily yields Theorem \ref{th main}. Finally,
section 7 gives a hypothesis on the random set $K$
that is weaker than (\ref{eq gosper-wiggly}), but still 
implies the conclusion of Theorem \ref{th any K}.
\section{Gosper Islands and Whitney tiles} \label{sec 2}
\begin{figure}
\centerline{ \psfig{figure=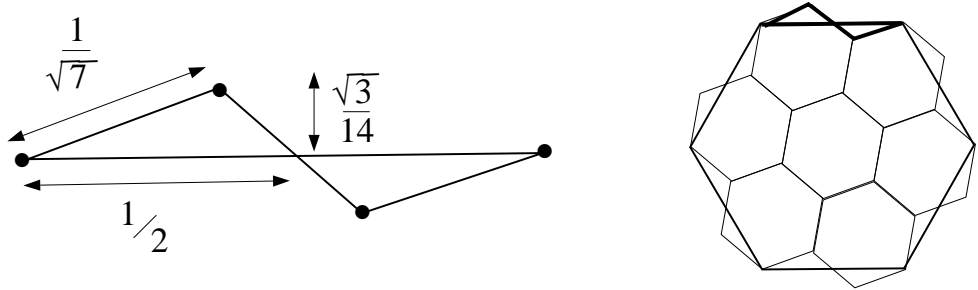,height=1.5in} }
 \caption{ \label{Sub} Substitution defining Gosper island }
\end{figure}

%
%
The standard hexagonal tiling of the plane is not self-similar,
but can be modified to obtain a self-similar tiling.
Replacing each hexagon by the union of seven smaller hexagons
(of area $1/7$ that of the original -- see Figure \ref{Sub}) yields a new tiling
of the plane by $18$-sided polygons; denote by $d_1$ the Hausdorff
distance between each of these polygons and the hexagon it approximates.
Applying the above operation to each of the seven smaller hexagons
yields a $54$-sided polygon with Hausdorff distance $7^{-1/2} \cdot d_1$ from 
the $18$-sided polygon,  which also has translates that tile 
the plane. Repeating this operation (properly scaled) ad infinitum,
we get a sequence of polygonal tilings of the plane, that
converge in the Hausdorff metric to a tiling of the plane by translates
of a compact connected 
set $G_0$ called the ``Gosper Island'' (see Gardner (1976)
and Mandelbrot (1982)).

\begin{figure}
\vbox{
 \centerline{ \hbox{
\psfig{figure=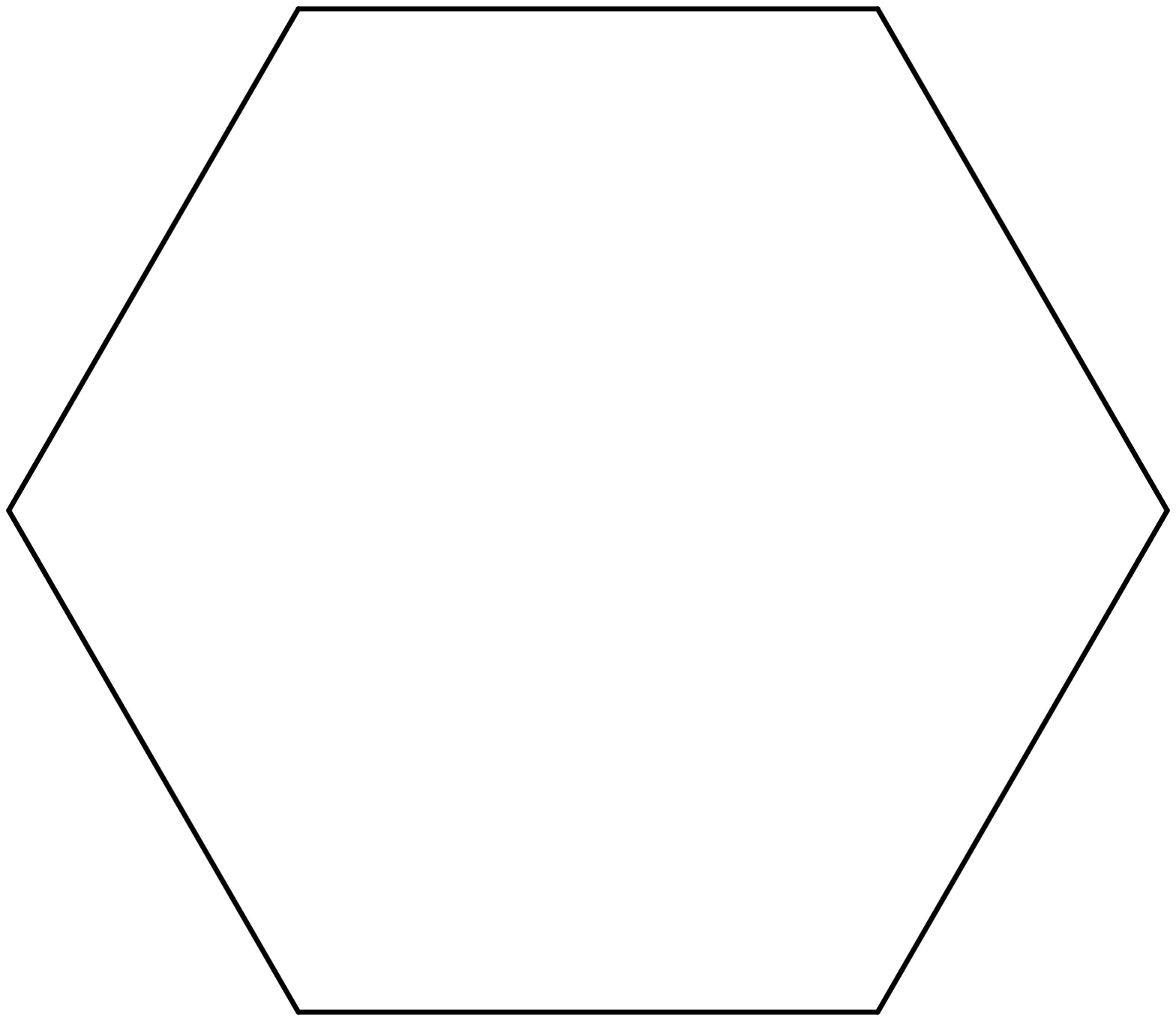,height=1.0in}
$\hphantom{xxx}$
\psfig{figure=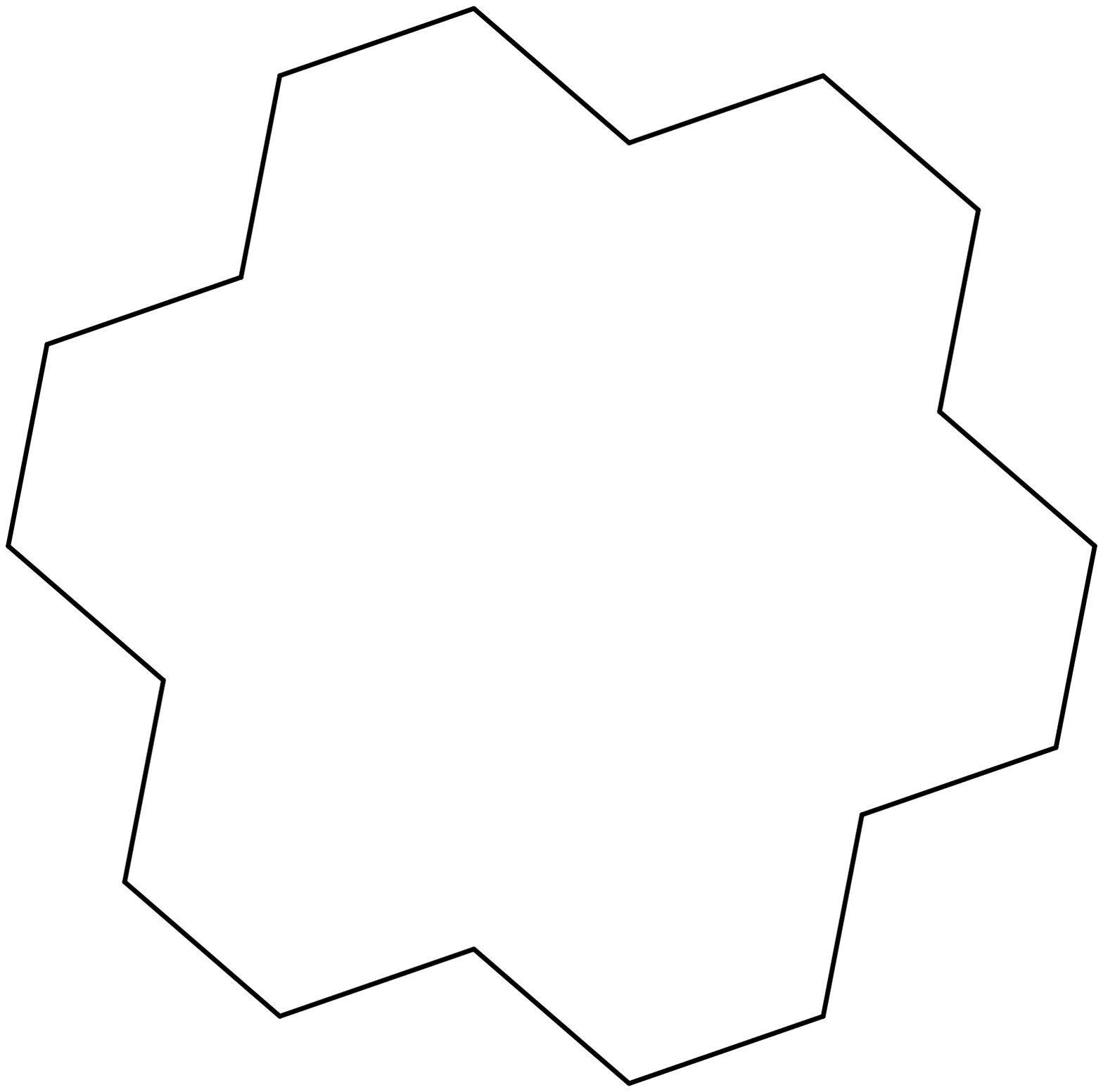,height=1.0in}
$\hphantom{xxx}$
\psfig{figure=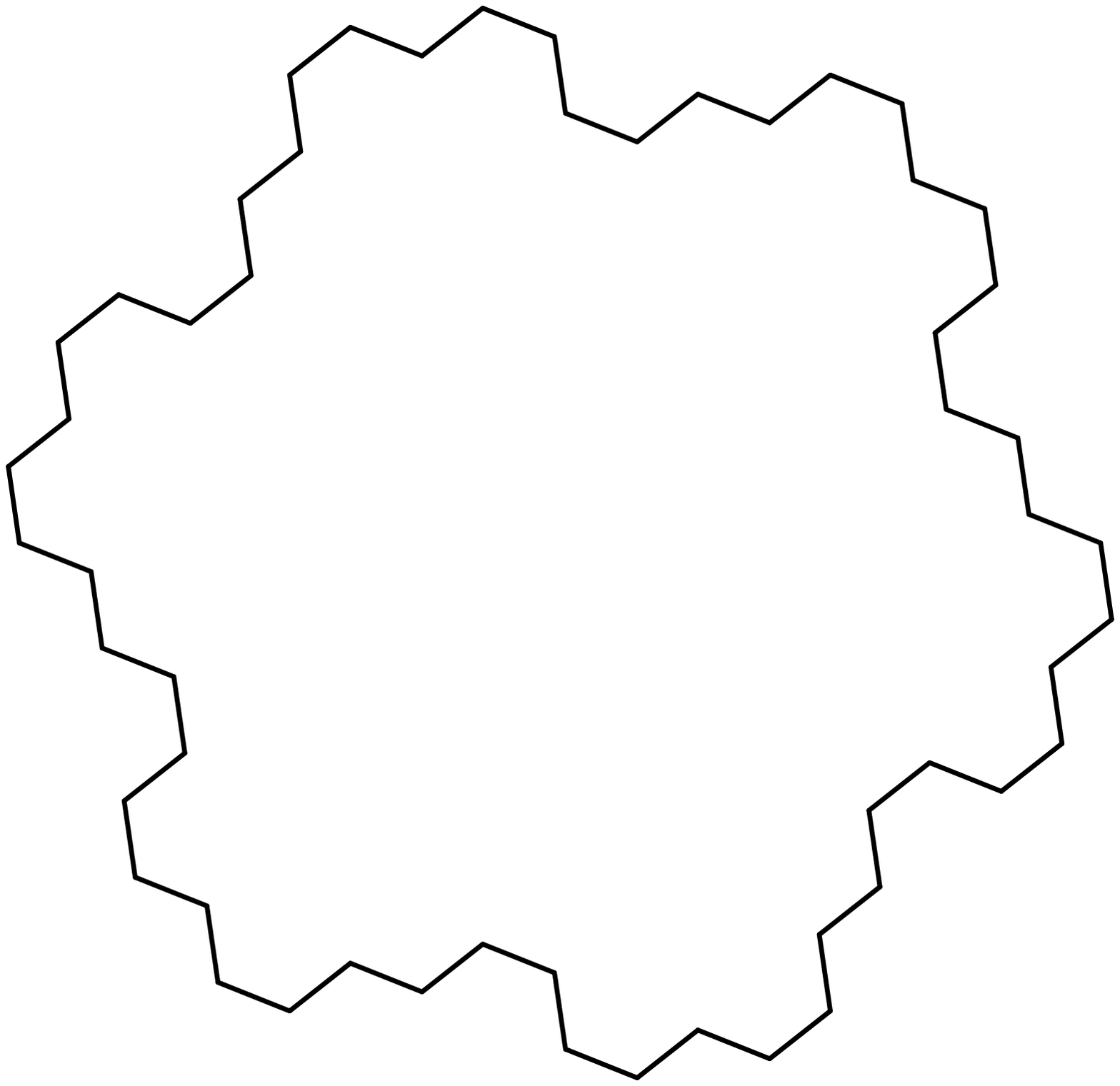,height=1.0in}
$\hphantom{xxx}$
\psfig{figure=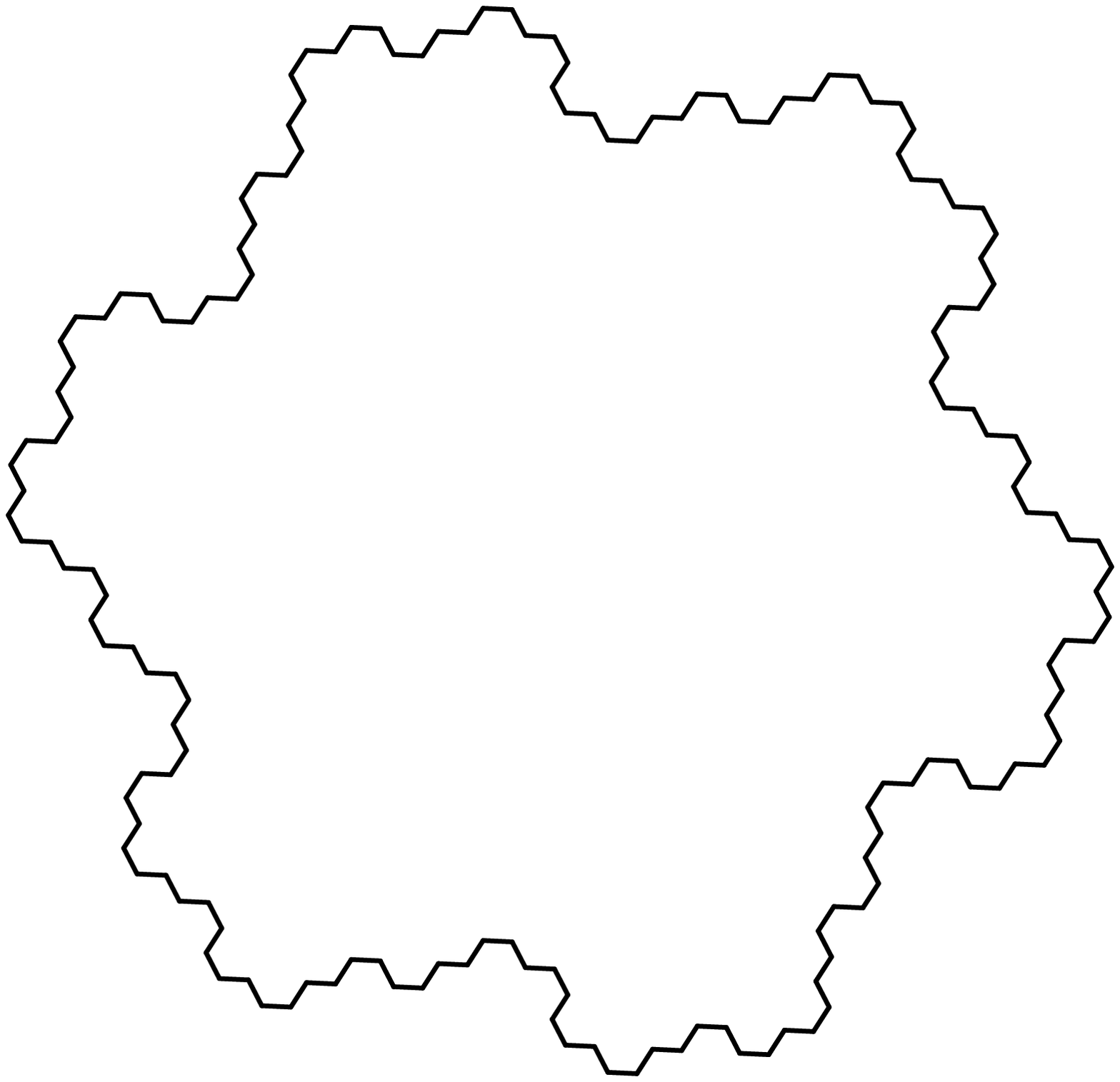,height=1.0in}
}}}
 \caption{ \label{Gosper} First four generation of the construction}
\end{figure}

\begin{figure}
\centerline{ \psfig{figure=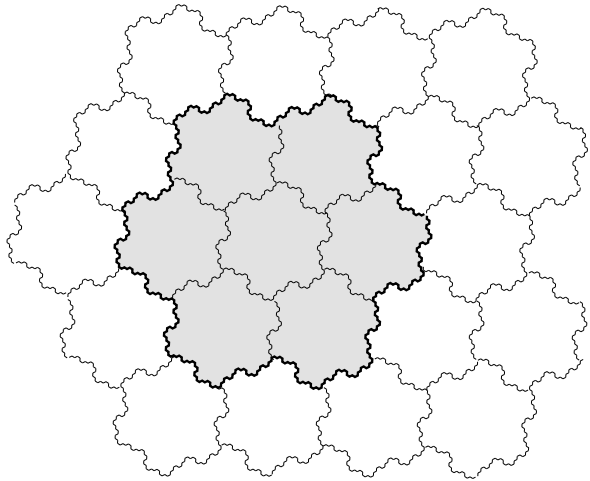,height=2.5in} }
 \caption{ \label{TileUnion} A self-similar tiling of the plane }
\end{figure}

\noindent{\bf Notation:} We normalize $G_0$ to be centered at the origin 
and have diameter 1.
Denote by $\D_0$ the set of translates of $G_0$ that form a tiling
of the plane (depicted in Figure \ref{TileUnion}).
This tiling is {\bf self-similar}, i.e., there is a complex number $\lambda$
with $|\lambda| >1$ such that  for each tile $G \in \D_0$,
the homothetic image $\lambda \cdot G$ is the union of tiles in $\D_0$.
(For the tiling by Gosper Islands, $|\lambda| = 7^{1/2}$.)
For each integer $n$, we denote by $\D_n$ the scaled tiling
$\{ \lambda^{-n} \cdot G \; : \: G \in D_0 \}$,
and let $\, \D = \cup_{n=0}^{\infty} \D_n$.
If $G \in \D_n$ we say that $G$ is a tile of index $n$ and write
 $||G||=n$. Every tile $ G \in \D_n$ is contained in a unique
tile of $\D_{n-1}$, denoted  $\, \pare(G)$. 
Each tile $G$  is centrally symmetric 
about a ``center point'' $z$;
for any $\theta >0$, denote by $\theta \odot G = z+ \theta \cdot (G-z)$
the expansion of $G$ by a factor $\theta$ around $z$.


\begin{figure}
\centerline{ \psfig{figure=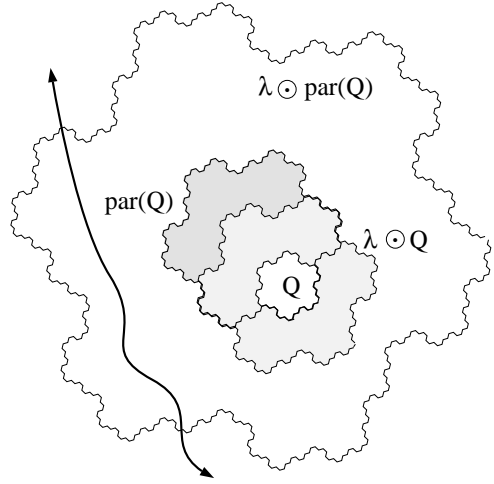,height=2.5in} }
 \caption{ \label{WhitTile} Boundary misses $\lambda \odot Q$, but 
 hits $\lambda \odot \pare(Q)$. }
\end{figure}

We record several simple properties of the tiling by Gosper Islands,
which will be useful later. 
\begin{enumerate} 

\item  \label{min-dist} 
 There is some minimal distance $d_0$ between
any two nonadjacent tiles of $\D_0$.

\item \label{eta}
 There is an $\eta_0 > 0$ such that any line segment 
of length $d_0$ must intersect $\core(G, 2\eta_0)$ for some $G \in \D_0$.
(The existence of $\eta_0$ follows by a compactness argument from the fact that
 $\partial G_0$ contains no straight line segments.)
\item \label{inrad}
The Gosper Island $G_0$ contains an open disk centered at the origin 
which in turn contains $\lambda^{-1} G_0$.

\item  \label{odot}
The blow-up $\lambda^3 \odot D$ contains $\lambda \odot 
\pare (D)$ for any $D \in \D$ (see  Figure \ref{WhitTile}).

\item  \label{neighbor}
 If $||G|| = ||G'|| - 1$ for neighboring tiles
$G$ and $G'$, then $\lambda \odot G$ contains $\lambda \odot G'$.
(See Figure~\ref{WhitTile}.)

\item \label{union}
 The blow-up $\lambda \odot G$ is contained in  $\bigcup (\lambda \odot G'')$
 where the union is over  all neighbors $G''$ of $G$
 of  index $||G||$.
\item \label{Jordan}
The boundary of $G_0$ is a Jordan curve.
To see this note that when we replace each segment or 
length $r$  by the 
three segments of the next generation, they remain 
within distance $r \sqrt{3}/14$ of the segment. Thus 
the limiting arc is within 
$$ r   \frac {\sqrt {3}}{14} \sum_{n=0}^\infty (\frac  1{\sqrt{7}})^n
 = {r \sqrt{21}\over 14(\sqrt{7} -1)} \approx r(0.198892),$$
of the segment. If $I_1, I_2, I_3$ are consecutive segments
of length $r$ then $\dist (I_1, I_3) = r$, so this shows the 
limiting arcs corresponding to them are at least distance
$r/2$ apart. Thus the boundary of the Gosper Island is a 
Jordan curve, indeed, is the image of the unit circle under 
a map $f$ satisfying
$$ \frac 1C \leq {|f(x) - f(y)| \over |x-y|^\alpha } \leq C,$$
where $\alpha = \frac 12 \log 7/\log 3$.
\item \label{annulus}
For any $\eta>0$, there is  a topological annulus 
with a rectifiable boundary, which separates $\core(G_0,\eta)$ from
 the boundary $\partial G_0$ of the Gosper island. \newline
(By the previous property,   the interior $G_0^\circ$ of $G_0$
 is simply connected, so this annulus can be obtained, for instance,
 by applying the Riemann mapping theorem.)
\end{enumerate}

\begin{figure} \centerline{ \psfig{figure=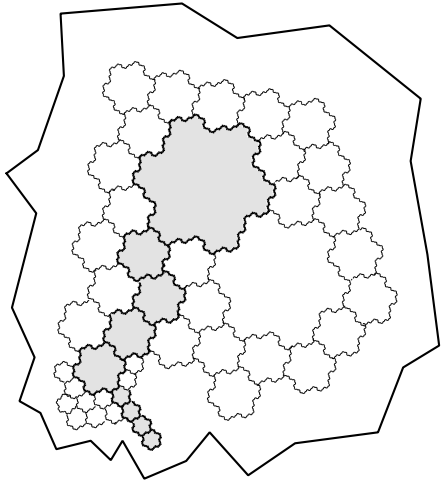,height=2.5in} }
 \caption{ \label{TileChain} Whitney decomposition and a chain of tiles }
\end{figure}

\noindent{\bf Definitions: }
Let $K$ be a compact connected subset of the plane. 
We  say that
$\, G \in \D \, $ is a {\bf Whitney~tile} for $\, K \,$ if
 $\, \lambda \odot G \,$ is disjoint from $\, K$,
 but $\, \lambda \odot \pare(G) \, $ intersects $\, K$.
(See Figure~\ref{WhitTile}.)
Let $W_K$ denote the set of Whitney tiles for $K$.
This collection is called a Whitney decomposition of
$K^c$, since it decomposes $K^c$ into a countable union 
of tiles (disjoint except for their boundaries) each with 
diameter comparable to its distance from $K$. See Figure
\ref{TileChain}.

A chain of adjacent tiles 
  $\{ G_1 , G_2 , \ldots , G_j \} \, $ in $\, W_K \,$ such that
  $\, G_i \subset 
\lambda^{5} \odot G_1 \, $ and   $\, ||G_{i}|| \geq ||G_{i-1}|| \,$ 
for all $\, i \in \{2, \ldots , j \} \, $ is called  a {\bf Whitney~chain}
(see Figure~\ref{TileChain}). 
Given $G \in W_K$, define $W_K^G \subset W_K$ to be the set of 
tiles $G'$ such that 
there is a Whitney chain 
  $\{ G_1 , G_2 , \ldots , G_j \}$  with $G_1 = G $ and $G_j = G'$.


Note the following property of the Whitney decomposition,
which holds for any connected component $\Omega$ of $K^c$: 
\begin{equation} \label{Whit-property}
 \mbox{ \rm 
For any open $V$ intersecting $\, \partial \Omega\, $, there
is a tile }  G_* \in W_K \mbox{ with } \lambda^{5} \odot G_*
\subset V. 
\end{equation}

\begin{lem} \label{lem 1}
If $G_1 , G_2 \in W_K$ are adjacent then $||G_1|| - ||G_2||$
is 0 or $\pm 1$.  
\end{lem}

\noindent{\sc Proof:} Suppose $\, ||G_1|| - ||G_2|| \geq 2$.  Let
$G$ be the tile of index $\, ||G_2||+1 \,$ that contains $\, G_1$,
and observe that $G$ is adjacent to $G_2$.  Then by Property \ref{neighbor}
of the Gosper tiling,  \newline
$\lambda \odot \pare (G_1) \subset
  \lambda \odot G \subset \lambda \odot G_2$ 
and maximality of $G_2$ is violated.  $\Cox$ 

\begin{lem} \label{lem 2}
Suppose $\C \subset W_K \cap \D_n$ is a collection of tiles whose
union topologically surrounds a smaller Whitney tile $G \in
W_K \cap \D_{n+k}$ where $k > 0$.  Then $\C$ surrounds a
a point of $K$.
\end{lem}

\noindent{\sc Proof:} The $k$-fold parent of $G$ is a tile $D \in\D_n$
which is surrounded by $\C$.  Applying Property \ref{union} inductively shows that 
the union of $\lambda \odot G'$ for $G' \in \C$ surrounds whatever part
of $\lambda \odot D$ it does not contain. Thus maximality of $G$
implies that $\lambda \odot D$ intersects $K$, and any point of 
intersection is  surrounded by $\C$.   $\Cox$

\begin{figure}
\centerline{ \psfig{figure=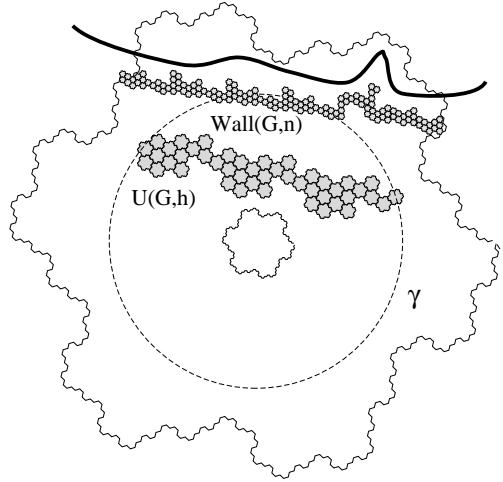,height=2.5in} }
 \caption{ \label{TileWall} The sets Wall(G,n) and U(G,h). }
\end{figure}

\begin{lem} \label{lem 3}
Suppose that $ \, G \in W_K \, $ and that there is a Whitney chain from 
 some larger tile outside $\, \lambda^{5} \odot G \, $ to $\, G$.  For any
$\, n > ||G|| \, $ define $\, \Wall (G,n) \,$ to be the set
 $ \, W_K^G \cap \D_n$.  
(See Figure~\ref{TileWall}.)  Let 
$$
E_n = \bigcup \{D \, : \,  D \in \Wall (G,n) \}
      \cup \partial (\lambda^{5} \odot G) \, .
$$
Then $ \, E_n \,$ is a connected set which topologically separates
 $\, G \,$ from
$\, K \,$.  Furthermore, If $\, \Gamma \, $ is a Jordan curve separating 
$\, G \, $ from the 
complement of $\, \lambda^{5} \odot G \, $, then every component of 
$\, \bigcup \{D \, : \,  D \in \Wall (G,n) \} \,$ intersecting the 
domain bounded by $\, \Gamma \,$ also intersects $\Gamma$.  
\end{lem}


\noindent{\sc Proof:} By Lemma \ref{lem 1} any path which connects 
$G$ to $K$ must hit Whitney tiles of every index 
larger than $||G||$. Thus any such path either hits 
$\Wall(G,n)$ or must leave $\lambda^{5} \odot G$, proving that 
$E_n$ separates $G$ from $K$.  Connectedness follows from the last
assertion of the lemma for $\Gamma = \partial (\lambda^{5} \odot G)$,
so it remains only to prove the last assertion.

Suppose to the contrary that there is a component $U$
of $\, \bigcup \{D \, : \,  D \in \Wall (G,n) \} \,$ which intersects
the domain bounded by $\Gamma$ but is disjoint from $\Gamma$ itself. 
 The union of all 
Whitney tiles which are in the unbounded component of 
$U^c$ and are adjacent to $U$ is a connected set.
By Lemma \ref{lem 1}  all of these tiles   have index  
$n-1$ or $n+1$. By connectedness and 
Lemma \ref{lem 1}, they must all have a single
index.  Suppose they all have index
$n+1$. Since tiles in $U$ can be connected to $G$ by 
Whitney chains which don't cross any tile of index $n+1$, this means
$G$ is in a bounded component of the complement of $U$, which 
contradicts our assumption that $G$ could be connected by a 
Whitney chain to a larger tile outside $\lambda^{5} \odot G$.
Thus the adjacent tiles must all have index $n-1$.
But then by Lemma \ref{lem 2}
these adjacent tiles must also surround a point of 
$K$, which implies that $K$ is not connected, another 
contradiction.   $\Cox$

The next two lemmas are needed in order to
show that if ``major portions'' of a wall of Whitney tiles
 can be covered by a thin
strip, then $K$ must intersect the core of an appropriate tile $G''$
 without $\eta$-surrounding it;
the latter event is controlled by the hypothesis of Theorem~\ref{th any K}. 
%
\begin{lem} \label{lem short segment}
Fix $G' \in \D$ and $\beta \in (0, \eta_0)$. Let 
$\hat U$ be any connected set intersecting both \newline
$\partial (\lambda^{5} \odot G')$ and $\lambda^3 \odot G'$.  Suppose that
$\hat U \cap (\lambda^{5} \odot G')$ is contained in an infinite open strip
of width $2 \beta \diam (G')$.  
Then there is a tile $G''$ contained in $\lambda^{5} \odot G'$ and 
of the same index as $G'$, such that 
 $\hat U$ intersects $\core(G'', 2\eta_0-2\beta)$.
\end{lem}

\noindent{\sc Proof:}  Pick a point $x \in \hat U \cap (\lambda^3 \odot G')$
and choose $y \in \hat U \cap \partial (\lambda^4 \odot G)$ connected
to $x$ inside $\hat U \cap (\lambda^4 \odot G)$.
By Property \ref{min-dist} of the tiling, the segment
$\overline{xy}$ has length at least $\, d_0 \cdot \diam (G')$,
so by Property \ref{eta} of the tiling, there is a tile
 $G''$ of the same index as $G'$, such that $\overline{xy}$
 contains some point $z \in  \, \core(G'', 2 \eta_0)$. 
By Property \ref{inrad} and convexity of disks, $z \in  \lambda^5 \odot G$,
and therefore $G'' \subset \lambda^5 \odot G$.
%
Observe that $\dist (z, \hat U) < {2 \beta \diam (G')}$, 
  for if not, removing the open disk centered at $z$ of radius 
$2 \beta \diam (G')$ from the infinite strip would 
contradict the connectedness of $\hat U$.   
This observation implies the assertion of the lemma.
    $\Cox$
 
Let $G$ be any Whitney tile with $\lambda^{5} \odot G$ not
containing all of $K$.  Let $\gamma$ be a circle 
centered at $\cen (G)$
which separates $\lambda^4 \odot G$ from
$\partial (\lambda^5 \odot G)$. (Such a circle exists
by Property \ref{inrad} of the tiling.)  For any
positive integer $h$, let $U(G,h)$ be the union of
all tiles $D \in W_K^G$ of index $||G|| + h$ such that $D$ intersects the
disk bounded by $\gamma$ (see Figure~\ref{TileWall}).  

\begin{lem} \label{lem wiggly -> beta}
Choose $a_2$ so that $\lambda^{3 - a_2} \leq \eta_0 / 2$.
With $G$ as above, let $G'$ be a tile in $W_K^G$ with 
$||G|| < ||G'|| < 
||G|| + h - a_2$  such that $\lambda^{5} \odot G'$ is contained in
the disk bounded by $\gamma$.
Suppose that  $U (G,h) \cap (\lambda^{5} \odot G')$ is covered by 
an open strip of width $2 \beta \diam G'$ with $\beta < \eta_0 /4$.
Then
there is a tile
 $G''  \subset \lambda^{5} \odot G'$  of the same index as $G'$,
 such that $K$  intersects $\core (G'', \eta_0) $ without 
$\eta_0$-surrounding $G''$.
\end{lem}

\noindent{\sc Proof:}  By Lemma~\ref{lem 3}, 
$U (G,h) \cup \gamma$ is connected
and therefore satisfies the hypotheses of Lemma~\ref{lem short segment};
let $G''$ be a tile as in the 
conclusion of that lemma.
Since $2\eta_0-2\beta > 3 \eta_0/2$, we can
pick a point $u$ in $U(G,h) \cap \core(G'', 3 {\eta_0}/2)$.
This clearly prevents $K$ from $\eta_0$-surrounding $G''$.
For any Whitney tile $D$ of index $||G||+h$, the blow-up 
$\lambda^3 \odot D$ intersects $K$ (by Property \ref{odot} of the tiling).
Since $U (G,h) \cap (\lambda^{5} \odot G')$ is a union of tiles
of index $||G||+h$, it follows that
$$
\dist(u,K) < |\lambda|^{3-||G||-h} \leq |\lambda|^{3-a_2-||G'||} \leq
  {\eta_0 \over 2} \diam (G'') \, ,
$$
by the choice of $a_2$.
Therefore $K$ intersects $\core(G'',\eta_0)$.
$\Cox$

\section{A tree of Whitney tiles} \label{sec 3}

Fix a compact, connected $K \subset \complex$, a tile $G_* \in W_K$,
and a positive integer, $h$.  We construct a tree $T = T(K , G_* , h)$ 
of Whitney tiles.  The root of $T$ is $G_*$ and the remaining
generations of $T$ are defined recursively as follows.  

Assume $T$ has been defined up to generation $n$ 
and for each $G$ in $T_n$, the $n^{th}$ generation of $T$, define 
$\widetilde T_{n+1} (G)$ 
to be the set of tiles $D$ with the following properties:
\begin{quote}
1.  $||D|| = \gs + (n+1)h$;

2.  $D \in W_K^G$;

3.  $\lambda^{5} \odot D \subset \lambda^{5} \odot G$.
\end{quote}
Let $ T_{n+1} (G)$ be 
a subcollection of $\widetilde T_{n+1} (G)$ which has maximal cardinality
 among all subcollections $\C$ for which
the expanded tiles $\{ \lambda^{6} \odot D : D \in \C \}$ are disjoint.
By maximality, $\bigcup\{ \lambda^{7} \odot D \, : \,  D \in T_{n+1} (G) \}$
contains all tiles in $\widetilde T_{n+1} (G)$ and therefore   
\be \label{eq:prune}
|\widetilde T_{n+1} (G)| \leq |\lambda|^{14} |T_{n+1} (G)| \, .
\ee
The children of $G$ in $T$ are defined to be the collection $ T_{n+1} (G)$.
Some trivial inductive observations are that $T_n \subset \D_{nh + \gs}$,
that each $G \in \D_n$ is connected to $G_*$ by a Whitney chain, and
that the sets $\lambda^{5} \odot G$ are disjoint as $G$ runs over any $T_n$.

\noindent{\bf Some tree terminology:} $\:$  Let $V$ be a countable set.
\begin{description}
\item{(i)} 
A mapping {\bf T} from a probability space ${\cal S}$ to the set
of trees on the vertex set $V$ is {\bf measurable} with respect to a 
$\sf $ $ \F$ on ${\cal S}$, if for any pair $\{v,v'\} \subset V$,
the event $\left[\{v,v'\} \mbox{\rm is an edge of {\bf T} } \right]$
is in $ \F$.
\item{(ii)}
  For any tree $T$ with vertex set contained in $V$,
and any element $v \in V$, define
$\trunc_v (T)$ to be null if $v$ is
not a vertex of $T$, and otherwise let $\trunc_v (T)$ be $T$ with the 
part below $v$ removed; more precisely, the vertices of $\trunc_v (T)$ are
the  vertices of $T$ not separated from the root by $v$, and the edges
are the edges of $T$ spanning pairs of vertices in this smaller vertex set.
\end{description}

For any $G \in \D$, let $\F_G$
denote the $\sf$ generated by the events
$\{ D \cap K \neq \emptyset \}$ for all tiles $D$ for which either
$||D|| \leq G$ or the interior of $D$ is disjoint from
$\lambda^{5} \odot G$. 

\begin{lem} \label{lem tree-dep}
On the event $||G|| = nh + \gs$, the random variable $\trunc_G \circ T$
is measurable with respect to $\F_G$.
\end{lem}

\noindent{\sc Proof:}  Suppose $||G|| = nh + \gs$ and consider an
event of the form 
$$\{ \{ D , D' \} \mbox{ is in the edge set of } \trunc_G (T) \} ,$$ 
where $D$ and $D'$ are tiles of index $mh + \gs$ and $(m+1) h + \gs$
respectively, \newline
 with $\lambda^{5} \odot D' \subset \lambda^{5} \odot D$.  
If $m \geq n$ and $\lambda^{5} \odot D$ is not disjoint from
$\lambda^{5} \odot G$, then the edge $\{ D , D' \}$ cannot be in 
$\trunc_G (T)$.  If $m < n$ or $\lambda^{5} \odot D$ is
disjoint from $\lambda^{5} \odot G$ then the event that
$\{ D , D' \}$ is an edge of $T$ is the union of events witnessed
by particular Whitney chains of tiles, all tiles being either
disjoint from $\lambda^{5} \odot G$ or of index at most 
$||G||$, so the event is measurable with respect to $\F_G$.   $\Cox$

The next lemma  requires the traveling salesman theorem
described in the next section, so its proof is delayed until
Section~\ref{sec 5}.
\begin{lem} \label{lem main}
Assume the random set $K$ satisfies the hypotheses of 
Theorem~\ref{th any K}.  Fix any  tile $G_*$ and $h > 0$
and let $T$ be the random tree $T(K , G_* , h)$.  
There are constants $\clol , \cllo > 0$ such that for any tile 
$G \in \D_{nh+\gs}$,
\begin{equation} \label{eq main}
\E [\#  T_{n+1} (G) \| \F_G ] \geq \clol h |\lambda|^h - \cllo |\lambda|^h 
\end{equation}
on the event that $G_* \in W_K$, the tile $G$ is in $T_n$,
 and $\lambda^{5} \odot G_*$ does not contain $K$.  The constants 
$\clol$ and $\cllo$ depend only on $\clll$.
\end{lem}

To prove Theorem \ref{th any K}, we also need two 
general lemmas concerning  trees.
 Define the 
{\bf boundary} $\partial T$ of the infinite rooted tree $T$ to be the set of
infinite self-avoiding paths from the root.  
The next lemma is implicit in Hawkes (1981) and can be found
in a stronger form in Lyons (1990). For convenience,
we include the short proof.

\begin{lem} \label{lem lyons}
Let $T$ be an infinite rooted tree. Given constants 
$C>0$ and $ \theta>1$, put a 
metric on $\partial T$ by 
\begin{equation} \label{eq: metric}
 \dist (\xi , \xi') = C  \theta^{-n} \mbox{ if }
   \xi \mbox{ \rm and } \xi' \mbox{ share exactly } n \mbox{ edges.}
\end{equation}
Suppose that independent percolation with parameter $p \in (0,1)$
is performed on $T$, i.e., each edge of $T$ is
 erased with probability $1-p$ and retained with probability $p$,
independently of all other edges.
If 
$$
 \dim(\partial T) < \alpha = {\log (1/p) \over \log \theta }
$$ 
then
with probability 1, all the connected components of retained edges in
$T$ are finite.
\end{lem}
\noindent{\sc Proof:}
It suffices to show that the connected component of the root
is finite almost surely. For any vertex $v$ of $T$,
denote by $|v|$ the number of edges between $v$ and the root.
By the dimension hypothesis
and the definition of the metric on $\partial T$,
there must exist cut-sets $\Pi$ in $T$ for which the
$\alpha$-dimensional cut-set sum
$$
\sum_{v \in \Pi} \theta^{-|v| \alpha} = \sum_{v \in \Pi} p^{|v|} 
$$
is arbitrarily small. But for any cutset  $\Pi$, the right-hand side is
the expected number of vertices in $\Pi$ which are connected to the root
 after percolation; this expectation bounds the probability
that the connected component of the root is infinite.
$\Cox$

The next lemma formalizes the notion of a random tree
which ``stochastically dominates'' the family tree of a 
branching process. We require the analogue of a filtration in our setting.

\noindent{\bf Definition:} $\:$ 
Let $V$ be a countable set and let $T$ be a random tree 
with vertex set contained in $V$, i.e., $T$ is a 
measurable mapping from some 
probability space $\, \langle {\cal S}, {\cal A}, \P \rangle$
to the set of trees on the vertex set $V$.
Say that $\sf$s $\{\F_v \, : \, v \in V\}$ on ${\cal S}$ 
form a {\bf tree-filtration}
if for any $v,w \in V$ and any $A \in \F_v$,
the event 
$A \cap \{ w \mbox{ \rm is a descendant of } v \mbox{ in } T \} \, $
    is  $\F_w$-measurable.

\begin{lem} \label{lem tree}
Let $V$ be a countable set and let $T$ be a random tree 
with vertex set contained in $V$. We assume that $T$ is rooted
at a fixed $v_* \in V$.   
Assume that  $b > 1$ and a tree-filtration $\{ \F_v \, : \,  v \in V \}$ 
 exists such that $\trunc_v (T)$ (defined before Lemma~\ref{lem tree-dep}) 
is $\F_v$-measurable for each $v \in V $, and the conditional expectation
$$
\E ( \mbox{number of children of } v \mbox{ in } T \| \F_v) \geq b .
$$
If every vertex of $T$ has at most $M$ children and at least
$m$ children, $m \geq 0$, then
\begin{enumerate}
\item The probability that $T$ is infinite is at least $1-q > 0$,
where $q$ is the unique fixed point in $[0,1)$ of the polynomial
$$
\psi (s) = s^m + {b-m \over M-m} (s^M - s^m) \, .
$$
  (Observe that $q = 0$ when $m > 0$.)
\item $\P (T \mbox{ is infinite } \| \tilde \F_{v_*} ) \geq 1 - q$ for any
$\tilde \F_{v_*} \subset \F_{v_*}$.
\item If $\partial T$ is endowed with the metric (\ref {eq: metric}), 
then $\dim (\partial T) \geq \log b / \log  \theta$ with probability
at least $1 - q$.
\end{enumerate}
\end{lem}

\noindent{\sc Proof:}
\noindent{\bf 1.}   Let $ |T_n|$ be the size of the $n^{th}$ generation
$T_n$ of $T$. 
 and let 
$\psi_n (s)$ denote the $n$-fold iterate of $\psi$.  
We claim that for $s \in [0,1]$,
\begin{equation} \label{eq claim 1}
\E s^{ |T_n|} \leq \psi_n (s) . 
\end{equation}
When $n = 1$, convexity of $x \mapsto s^x$ implies that
$$
s^{|T_1|} \leq s^m + 
{|T_1| - m \over M-m} (s^M - s^m)
$$
 and the claim follows by 
taking expectations:
$$\E s^{|T_1|} \leq s^m + {\E |T_1| - m \over M-m} (s^M - s^m)  
   \leq \psi_1 (s)$$
since $s^M - s^m \leq 0$.  For $n > 1$ proceed by induction.  Let 
$|T_{n+1} (v)|$ be the number of children of $v$ if $v \in T_n$ and zero 
otherwise, and use the argument from the $n=1$ case to 
see that $\E (s^{ |T_{n+1} (v)|} \| \F_v ) \leq \psi (s)$ on the event
$v \in T_n$ (which is an event in $\F_v$).  Giving $V$ an arbitrary
linear order (denoted ``$<$''), we  have in particular
$$\E (s^{ |T_{n+1}(v)|} \| T_n \mbox{\rm and }
    T_{n+1} (w) \mbox { for } w < v \mbox { in } T_n) 
        \leq  \psi (s) \, .
$$
for $v \in T_n$. Since 
$$
s^{ |T_{n+1}|} =\prod_{v \in T_n} s^{ |T_{n+1}(v)|} \, , \;
\mbox{ this  yields }
 \; \; \E (s^{ |T_{n+1}|} \| T_n)   \leq  \psi (s)^{ |T_n|} \, .
$$
Taking expectations and applying the induction hypothesis with
$\psi (s)$ in place of $s$ gives
$$
\E s^{ |T_{n+1}|} \leq  \E \Big(\psi (s)^{ |T_n|}\Big) 
  \leq \psi_n (\psi (s)) = \psi_{n+1} (s) \, ,
 $$
proving the claim (\ref{eq claim 1}).  

                         From (\ref{eq claim 1}) we see that 
$\P ( |T_n| = 0) \leq E q^{ |T_n|} \leq \psi(q) =q $, 
establishing the first conclusion of the lemma. 

\noindent{\bf 2.}
By copying  the derivation of~(\ref{eq claim 1}), inserting an extra
 conditioning on $\tilde \F_{v_*}$, one easily verifies
 that $ \E (s^{ |T_n|} \| \tilde \F_{v_*}) \leq \psi_n (s)$, and 
the rest of the argument is the same as in the first part.   

\noindent{\bf 3.}
Let $T' (v)$ be the connected component of the 
subtree of $T$ below $v$ after removing each vertex of $T$ below $v$
independently with probability $1-p$.  For $p > 1/b$, let $q_p
\in (0,1)$ solve $q_p = 1 + (bp/M) (q_p^M - 1)$.  
We apply the second part of the lemma to $T' (v)$ conditioned
on $\F_v$ to see that
$$\P (T' (v) \mbox{ is infinite } \| \F_v ) \geq 1 - q_p$$
for $v \in T$.  By Lemma \ref{lem lyons}, the
event $\{ \dim (\partial T) < |\log p| / \log  \theta \}$
is contained up to null sets in the event $\{ T' (v) 
\mbox{ is finite for all } v \in  T_n \}$.  Thus 
\begin{eqnarray*}
\P \Big( \dim (\partial T) < |\log p| / \log  \theta \cond  T_n \Big) 
& \leq & 
   \P \Big( \cap_{v \in  T_n} T' (v) \mbox{ finite } \cond T_n \Big) \\[2ex]
& = & \prod_{v \in T_n} \P \Big( T' (v) \mbox{ finite } \cond  T_n  \, , \, 
   T' (w) \mbox{ finite for all } w < v \mbox{ in } T_n \Big) \\[2ex]
& \leq & q_p^{ |T_n|}
\end{eqnarray*}
since each event conditioned on is in the corresponding $\F_v$.
Taking expectations yields
$$\P \Big( \dim (\partial T) < {|\log p| over \log  \theta} \Big)
 \leq \psi_n (q_p) .$$
Since $q_p < 1$ for each $p > 1/b$, this goes to $q$ as
$n \rightarrow \infty$, proving the last conclusion 
of the lemma. 
$\Cox$ 

\noindent{\sc Proof of Theorem}~\ref{th any K}: 
 Put a metric on $\partial T$ by 
$$ \dist (\xi , \xi') = |\lambda|^{-(n+1)h - \gs} \mbox{ if }
   \xi \mbox{\rm and } \xi' \mbox{ share exactly } n \mbox{ edges.}
$$
Each $\xi = (G_1 , G_2 , \ldots) \in \partial T$ defines a 
unique limiting point $\phi (\xi) \in \frK$ which is the 
decreasing limit of the set $\lambda^{5} \odot G_n$.  If
$\xi = (G_1 , G_2 , \ldots)$ and $\xi' = (G_1' , G_2' , \ldots)$ 
share exactly $n$ edges, then by definition of $ T_{n+1}$, the expanded tiles
$\lambda^{6} \odot G_{n+1}$ and $\lambda^{6} \odot G_{n+1}'$ 
are disjoint.  Since $\phi (\xi) \in \lambda^{5} \odot G_{n+1}$
and $\phi (\xi') \in \lambda^{5} \odot G_{n+1}'$, it 
follows from Property \ref{min-dist} of the tiling that 
$$|\phi (\xi) - \phi (\xi')| \geq d_0 |\lambda|^{-(n+1)h - \gs} .$$
Thus 
$$|\phi (\xi) - \phi (\xi')| \geq d_0 \cdot \dist(\xi , \xi')$$
and since the range of $\phi$ is included in $\frK \cap
\lambda^{5} \odot G_*$ it follows that
\begin{equation} \label{eq local dim}
\dim (\frK \cap \lambda^{5} \odot G_*) \geq \dim (\partial T) .
\end{equation}
    From  Lemma \ref{lem tree} and the conclusion of Lemma~\ref{lem main}, we
see that 
\begin{equation} \label{eq dim bound}
\dim (\partial T (K , G_* , h)) \geq { \log ( \clol h |\lambda|^h - 
   \cllo |\lambda|^h ) \over h \log |\lambda|} 
\end{equation}
with probability 1, on the event that $G_* \in W_K$ and $\lambda^{5}
\odot G_*$ does not contain $K$.  Choose $h$ 
to maximize the RHS of~(\ref{eq dim bound}).  Since the maximum 
is greater than 1, there is an $\eps > 0$ for which 
$$\dim (\partial T (K , G_* , h)) \geq 1 + \eps$$
with probability 1 on this event. Finally, let $\Omega$ be any connected 
component of $K^c$.
 By  property (\ref{Whit-property}
 of the Whitney
decomposition, 
for any open $V$ intersecting $\partial \Omega$,
 there is a tile $G_* \in W_K$
with $\lambda^{5} \odot G_* \subset V$, and the theorem
follows from~(\ref{eq local dim}). 
  $\Cox$

\noindent{\bf Remark:} A planar domain $\Omega$ is called
a {\bf John domain} if
there is a base point $z_0 \in \Omega$
and a constant $C >0$  so 
that any point $x \in \Omega$ can be joined to $z_0$ by a 
curve $\gamma_x \subset \Omega$ so that 
$ \dist(z, \partial \Omega) \geq C|x-z| $ for any $z \in \gamma_x$. 
John domain were introduced by Fritz John in 1961, and
some basic facts about them can be found  
in N{\"a}kki and V{\"a}is{\"a}l{\"a} (1994).

With the notation of the above
proof, if $G_* \in W_K$ is contained in a component $\Omega$
of $K^c$, choose for every tile  $G \neq G_*$ in the tree
$T (K , G_* , h)$, a  Whitney chain leading to $G$ from its unique 
ancestor in the previous generation of the tree.
For each tile $G'$ in this chain, there is an open disk containing
it which is contained in $\lambda \odot G'$ (by Property \ref{inrad}
of the tiling).
The union of all these open disks as $G'$ runs over the chosen Whitney
chain for $G$ and $G$ runs over $T (K , G_* , h)$,
 is a John domain $\Omega_{\rm J}$ satisfying
$\dim (\partial \Omega \cap \partial \Omega_{\rm J}) \geq  1 + \eps$. 

\section{The traveling salesman theorem} \label{sec 4}

Given a set $E$ in the plane and another bounded
plane set $S$, we define
$$ \beta_E (S) = (\diam (S))^{-1} \inf_{ L \in \cal L}
   \sup_{z \in E \cap S} \dist (z, L),$$
where $\cal L$ is the set of all lines  $L$  intersecting $S$.
\begin{th}[Jones 1990] \label{th jones}
If $E \subset \complex $ then the length of the shortest
connected curve $\Gamma$ containing $E$ is bounded between
(universal) constant multiples of
$$\diam(E) + \sum_Q \beta_E (3 \odot Q)^2 \diam (Q) , $$
where the sum is over all dyadic squares in the plane and
$3 \odot Q$ is the union of a 3 by 3 grid of congruent squares
with $Q$ as the central square.
\end{th}

A simpler proof of this Theorem, and an extension to higher dimensions,
are given in Okikiolu (1992).
The theorem easily implies that the 
length $|\Gamma|$ of any curve $\Gamma$ which passes 
within $r$ of every point of $E$ satisfies
\begin{equation} \label{eq jones}
\diam(E) + \sum_{\diam(Q) \geq r} \beta_E 
   (3 \odot Q)^2 \diam (Q) \leq \clo |\Gamma| ,
\end{equation}
where the sum is over all dyadic squares in the plane with
diameter at least $r$.
For every set $S$, there is a dyadic square $Q$ of side length at most
$2 \diam (S)$ for which $S \subset 3 \odot Q$.  
Picking $S = \lambda^{5} \odot G$ for some tile $G$ and
$Q$ accordingly, we get
$$\beta_E (\lambda^{5} \odot G)^2 \diam (G) \leq 
   9 |\lambda|^{-10} \beta_E (3 \odot Q)^2 \diam (Q) $$
and since each expanded square $3 \odot Q$ contains a bounded number of
expanded tiles $\lambda^{5}
\odot G$ for tiles $G$ with $\sqrt{2} |\lambda|^{5} \diam (G) \geq 
\diam (Q)$, it follows that the length of any curve passing
within $r$ of every point of $E$ satisfies
\begin{equation} \label{eq jones2}
|\Gamma| \geq \cl \left ( \diam(E) + \sum_{\diam(G) \geq r} 
   \beta_E (\lambda^{5} \odot G)^2 \diam (G) \right ) .
\end{equation}
We require the following corollary, which uses an idea from
Bishop and Jones (1994).

\begin{cor} \label{cor BJ}
Let $\gamma$ be a Jordan curve with length denoted $|\gamma|$ and
let $\C$ be a collection of Whitney tiles of index $n$.  Let $U$
denote $\bigcup_{D \in \C} D$ and suppose that $\gamma \cup
U$ is connected.  Then there is a constant $\cll$
such that the cardinality of $\C$ is at least
$$\cll |\lambda|^n \left ( - |\gamma| + \sum_{G' \in \Xi (\C)} \beta_U
   (\lambda^{5} \odot G')^2 \diam (G') \right ) \; , $$
where $\Xi (\C)$ is the collection of tiles $G'$ of index at most $n$
for which $\lambda^{5} \odot G'$ intersects $U$.
\end{cor}

\noindent{\sc Proof:} Let $\C_\circ$ be the collection of circles of
radius $r := |\lambda|^{-n} $ centered at points $\cen (D)$
for $D \in \C$.  Since neighboring tiles in $\C$ give rise to
intersecting circles in $\C_{\circ}$, we see that $\Gamma := \gamma \cup
\bigcup_{\Theta \in \C_{\circ}} \Theta $ is connected and passes within
$r$ of every point of $\gamma \cup U$.  Furthermore, any
connected finite union of closed  curves is a closed
curve, and hence $\Gamma$ is a curve of length at most
$|\gamma| + 2 \pi \# (\C) |\lambda|^{-n} $.
Combining this with~(\ref{eq jones2}) shows that
$$\# (\C) \geq {|\lambda|^n \over 2 \pi } \, \Big (
   \cl \sum_{\diam (G') \geq r} \beta_U (\lambda^{5} \odot G')^2
   \diam (G')  - |\gamma| \Big ) .$$
Since all tiles in $\Xi (\C)$ have diameter at least $r$, this proves
the lemma.   $\Cox$

\section{ Expected offspring in the Whitney tree} \label{sec 5}

\noindent{\sc Proof of Lemma}~\ref{lem main}: 
  Fix $G_*$ and $h$ 
 as in the statement of the lemma and let $G$ be
any tile in $T_n$.  Let $\gamma$ be the circle separating
$|\lambda|^4 \odot G$ from $\partial (|\lambda|^5 \odot G)$,
which was used in Lemma \ref{lem wiggly -> beta}.
Let $\C$ be the 
collection of tiles $D \in W_K^G$ of index $||G|| + h$ intersecting the
disk bounded by $\gamma$.  The union of all tiles in
$\C$ is  the set $U= U (G,h)$
defined before Lemma \ref{lem wiggly -> beta}.
  We want to 
show that the expected cardinality of $ T_{n+1} (G)$ is large.
 Since the cardinality of $ T_{n+1}(G)$  is at least
$|\lambda|^{-14}$ times the cardinality of $\widetilde T_{n+1}(G)$
by (\ref{eq:prune}), and $\widetilde T_{n+1} (G) $ is a superset of $\C$,
 it suffices to show that
$$\E (\# \C \| \F_G) \geq \clol' h |\lambda|^h - \cllo' |\lambda|^h \, .$$
To do this, we will apply Corollary~\ref{cor BJ} to $\C$,
 so that the set $U$ defined in that corollary 
is the same as $U(G,h)$ defined above.  

\begin{figure}
\centerline{ \psfig{figure=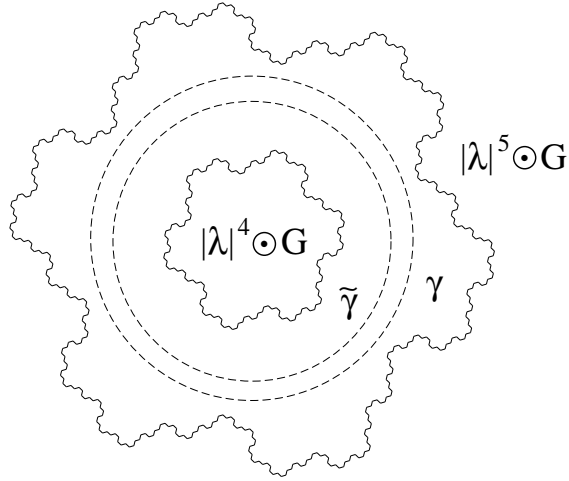,height=2.5in} }
 \caption{ \label{GammaCir} The circles $\gamma$ and $\tilde \gamma$ }
\end{figure}

We will be able to bound from below the
summands in Corollary~\ref{cor BJ} for most, but not all,
 ``intermediate-sized'' tiles $G'$.
Pick an integer $a_3>1$ so that $|\lambda|^{3-a_3} < d_0$,
where $d_0$ is the minimal distance between nonadjacent tiles in $\D_0$.
Let $\wtg$ be a circle concentric with $\gamma$,
with a smaller radius:
$\rad(\wtg)=\rad(\gamma)-|\lambda|^{5-a_3}$ (see Figure \ref{GammaCir}).

For $a_3< j < h $, let $W_K^G (j)$ be the set of tiles 
$G' \in W_K^G$
 such that $ ||G'||=||G|| + j$ and $\lambda^5 \odot G'$ intersects
  the disk bounded by $\wtg$ . For any such tile $G'$
the blow-up $\lambda^5 \odot G'$ is contained in the disk bounded
by $\gamma$.

 Fix any tile $G' \in W_K^G(j)$ with $a_3 < j \leq ||G|| + h - a_2$,
where $a_2$ was specified in Lemma~\ref{lem wiggly -> beta}.
 That lemma implies 
\begin{eqnarray}
&& \P \left ( \beta_U (\lambda^{5} \odot G') \geq 
 {\eta_0 \over 4 } \cond \F_G \mbox {\rm and } W_K^G(j) \right ) 
                                      \label{eq asdf} \\[2ex]
& \geq & \P \Big[ \bigcap_{G''}
  \Big(  K \, \mbox{ \rm  $\eta_0$-surrounds } G''
  \, \mbox { \rm \underbar{or} } \; 
   K \cap \! \core(G'',\eta_0) = \emptyset \Big) 
   \cond \F_G \mbox {\rm and } W_K^G(j) \Big] \, \nonumber ,
\end{eqnarray}
where the intersection is over all tiles 
$G'' \subset \lambda^{5} \odot G'$
such that $||G''|| = ||G'||$.
The set of such tiles $G''$ for a fixed $G'$ has cardinality
$|\lambda|^{10}$. 
 Enumerating these and multiplying
conditional probabilities using the hypotheses of 
Lemma~\ref{lem main} (since the $\sigma$-fields $\F_G$  and  
  $\sigma(W_K^G(j))$ are contained in $\sigma(K \setminus G''^{\circ})$ )
 gives a lower bound of
$\clll^{|\lambda|^{10}} $
for~(\ref{eq asdf}), and implies that
$$
\E \Big( \beta_U (\lambda^{5} \odot G')^2 \cond
   \F_G \mbox{\rm and } W_G^K (j) \Big)
   \geq \left ( {\eta_0 \over 4 } \right )^2 
   \clll^{|\lambda|^{10}} = \cbeta>0 \, .
$$
(This is the definition of $\cbeta$).
  Since $\gamma$ is outside $\lambda^4
\odot G$, the distance from $\gamma$ to $\lambda^3 \odot G$
is at least $|\lambda|^3 d_0 \cdot \diam(G)$,
 by Property \ref{min-dist} of the tiling.
Therefore the distance from
$\wtg$ to $ \lambda^3 \odot G$ is at least
 $ (|\lambda|^3 d_0 - \lambda^{5-a_3} )\cdot \diam(G)$,
which is greater than  $ d_0 \cdot \diam(G) \,$ by the choice of $a_3$.
 By Lemma~\ref{lem 3}, 
the union of $\wtg$ with all the tiles in $W_K^G (j)$ is a connected
set. Since it intersects both $\lambda^3 \odot G$ and $\wtg$,
 it follows that  the cardinality of
$W_K^G (j)$ is  at least $d_0 |\lambda|^j$.  
Thus for each  $j \in (a_3, h-a_2]$ we have
$$\E \Big( \sum_{G' \in W_K^G (j)} \beta_U (\lambda^{5} \odot G')^2 \cond
   \F_G \mbox{\rm and } W_G^K (j) \Big)
   \geq \cbeta d_0 |\lambda|^j  .
$$
By Corollary~\ref{cor BJ}, 
$$\E (\# \C \| \F_G) \geq \cll |\lambda|^{(n+1)h+\gs} \left (
   -|\gamma| + \sum_{j = a_3+1}^{h - a_2}
    |\lambda|^{-nh-\gs - j} \cbeta d_0  |\lambda|^j \right ) .
$$
Summing gives 
$$\E ( \# \C \| \F_G) \geq \cll \left (
    \cbeta d_0 (h - a_2-a_3) 
   |\lambda|^h - |\gamma| \cdot |\lambda| ^{(n+1)h + \gs} \right )$$ 
which proves the lemma since $|\gamma| = 2 \pi |\lambda|^{4 - nh - \gs}$.
  $\Cox$ 

\section{The Brownian frontier: Proof of Theorem 1.1} 
\label{sec 6}
Let  $\P_x$ denote the law of a planar
Brownian motion $\{B(t) \}_{t \geq 0}$
started at $x$. We use $\P_0$ unless indicated explicitly otherwise.
Let $\tau_{\exp}$ be a positive random variable,
 independent of the Brownian
motion, which is exponential  of mean 1 (i.e., its density is $e^{-t}$).
  We will verify~(\ref{eq gosper-wiggly}) for $K = B[0,\tau_{\exp}]$.
by Brownian scaling, this will imply the first assertion of
Theorem~\ref{th main}.  

\noindent{\bf Notation: } for any compact planar  set
$S$, denote by $\tau_S = \min \{ t \geq 0 \, : \, B(t) \in S \}$
the first hitting time of $S$,
which is almost surely finite if $S$ has positive logarithmic capacity.
Given $\eta \in (0,1/10)$,
let $J_0$ be a rectifiable closed Jordan curve,
which is the exterior boundary of a topological annulus separating
$\core(G_0, \eta)$ from $\partial G_0$. (Here the constant $1/10$ can be 
replaced by any constant smaller than the inradius of $G_0$,
and the existence of $J_0$ is guaranteed by 
Property \ref{annulus} of the tiling.)
For the rest of this section, consider
 a homothetic image $G=z_{\rm cen}+r G_0$ of the Gosper Island $G_0$,
 with $r \in (0, 1)$ and $z_{\rm cen}$ in the plane.
Also,  denote by $J=z_{\rm cen}+r J_0$
the image of $J_0$ in $G$.
We must obtain estimates which are uniform in the location and scale of
$G$, as well as in the structure of the Brownian range outside $G$.
\begin{lem} \label{lem:fastloop}
For every $x \in \corG$,
\be \label{eq:loop}
\P_x \Big( B[0,\tau_J] \; \,
  \mbox{\rm $\eta$-surrounds } G \Big) \geq \cloop(\eta) >0 \, .
\ee
Furthermore, there exists $\cfast(\eta)>0$ such that
\be \label{eq:fast}
\P_x \Big( B[0,\tau_J] \;  \mbox{ \rm $\eta$-surrounds } G  \: \mbox{ \rm and }
\tau_J < \tau_{\exp} < \tauG \Big) \geq 
  \cfast(\eta) \P_x(\tau_{\exp} < \tauG)     \, .
\ee
\end{lem}

\noindent{\sc Proof:} The first estimate is immediate for $G_0$,
and the general case follows by  scaling.
For the second, observe that by Brownian scaling,  
$\; \inf_{y \in J} \P_y (\tauG > \diam(G)^2) \: $ is a 
positive constant depending
only on $J_0$, hence only on $\eta$. Also, clearly 
$\P(\tau_J < 1) >1/2$ and therefore
$$
 \P(\tau_J < \tau_{\exp} < \tau_J+ \diam(G)^2) > 
{e^{-2} \over 2} \diam(G)^2 \,   ,  \; \mbox{   since } \, \diam(G) <1. 
$$
 Applying (\ref{eq:loop}), lack of memory of exponential variables, and
the strong Markov property 
of Brownian motion at the stopping time $\tau_J$,
 then shows that the left-hand side of
(\ref{eq:fast}) is at least a constant multiple of $\diam (G)^2$.
On the other hand, for any $x \in G$ we have 
$\P_x(\tau_{\exp} < \tauG) \leq \E_x ( \tauG) \leq \diam(G)^2$.
This completes the proof.
$\Cox$

\begin{lem} \label{lem 3 parts}
 There exists $\clool (\eta) > 0$ such
that for any tile $G$, for any $x \in \corG$ and any $A \subset 
\partial G$,
$$
\P_x \Big(B[0,\tauG] \; \,
   \mbox{ \rm $\eta$-surrounds } G \; \mbox{\rm and } B(\tauG) \in A  \Big)
 \geq \clool(\eta) \P_x (B(\tauG) \in A) .
$$
\end{lem}
\noindent{\sc Proof:}  
 Recall that $J$ is a Jordan curve of finite length
separating $\corG$ from $\partial G$.  The Harnack principle (see,
e.g.,~Bass (1995, Theorem~1.20)) implies
that there is a constant \newline $\clolo = \clolo (J_0)$
 such that for any $y,z \in J$
and for any $A \subset  \partial G$, 
\be \label{eq:harnack}
\P_y (B(\tauG) \in A) \geq \clolo \P_z (B(\tauG) \in A) .
\ee
Therefore for any $x \in \corG$,
\begin{eqnarray}
&\P_x ( B[0,\tau_J] \; \, \mbox{ \rm $\eta$-surrounds } G \; \mbox{\rm and } 
 B(\tauG) \in A )  &  \nonumber \\[3ex]
&= \E_x \Big( \one_{\{ B[0,\tau_J] \; \mbox{ \rm $\eta$-surrounds } G \}}
 \cdot \P_{\! B(\tau_J)} [B(\tauG) \in A ] \Big) &     \label{eq:one-ave}
\end{eqnarray}
Applying the Harnack inequality (\ref{eq:harnack}) with $y=B(\tau_J)$
and then invoking the estimate
(\ref{eq:loop}) from the previous lemma, we find that the expression
(\ref{eq:one-ave}) is at least $\cloop(\eta) \clolo \P_z (B(\tauG) \in A)$,
for any $z \in J$. Finally, taking $z=B(\tau_J)$ and
 averaging with respect to 
$\P_x$ using the strong Markov property gives
$$
\P_x ( B[0,\tau_J] \; \, \mbox{ \rm $\eta$-surrounds } G \mbox{\rm and } 
 B(\tauG) \in A ) \geq \cloop(\eta) \clolo \P_x (B(\tauG) \in A) \, ,
$$
for any Borel set $A \subset \partial G$.
This proves the lemma.
$\Cox$

\begin{figure}
\centerline{ \psfig{figure=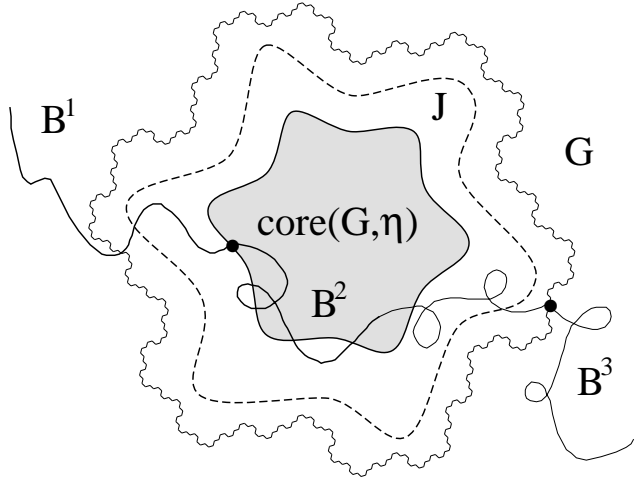,height=2.5in} }
 \caption{ \label{B1B2B3} The partition of the Brownian trajectory }
\end{figure}

Given $\eta>0$, we abbreviate $\tauc= \tau_{\core(G, \eta)}$ and
partition the Brownian trajectory into three
pieces:
\begin{enumerate}
\item Until the first time $\tauc$ that the path visits $\core(G, \eta)$.
 \newline
      Formally, define $B^{(1)}(t) = B(t \wedge \tauc)$ for $t \geq 0$,
  where $t \wedge s$ is shorthand for $\min \{t,s\}$.
\item From time $\tauc$ until the next visit to $\partial G$, denoted
 $\taucG=
\min \{ t \geq \tau_c \, : \, B(t) \in \partial G \}$. \newline
Define $B^{(2)}(t) = B((t+\tauc) \wedge \taucG)$ for $t \geq 0$.
\item After time $\taucG$. \newline 
Denote $B^{(3)}(t)=B(t+\taucG)$ for $t \geq 0$.
\end{enumerate}
The idea now is that $B^{(1)}$ and $B^{(3)}$ determine the Brownian range
outside $G^\circ$, and $B^{(2)}$ has a substantial chance
of $\eta$-surrounding $G$, even when we condition on its endpoints.
However, we still have to take the exponential killing into account.
Define the random variable
$$
I = \left\{ \begin{array}{lcr}
                    1 & \mbox{\rm if} & \tau_{\exp} < \tauc \\
                    2 & \mbox{\rm if} & \tauc \leq \tau_{\exp} < \taucG \\
                    3 & \mbox{\rm if} & \taucG \leq \tau_{\exp}   
\end{array} \right. \, 
$$
that indicates in which part of the motion the exponential killing occurred.
Finally, define
$$ 
\texp = \left\{ \begin{array}{lcr}
           \tau_{\exp} & \mbox{\rm if} & I=1 \\
           \tauc & \mbox{\rm if} & I=2 \\
           \tau_{\exp}-\taucG & \mbox{\rm if} & I=3 
\end{array} \right. 
$$   

\begin{pr} \label{pr:brown}
 For any $\eta \in (0, 1/10)$ there is a constant $\clll= \clll(\eta) > 0$
such that for  all homothetic
images $G=z_{\rm cen} +r G_0$ of $G_0$ with $r \in (0, 1)$
and $z_{\rm cen}$ in the plane:  
\begin{eqnarray} \label{eq gosper-pf}
   \P \Big ( B[0, \tau_{\exp}] \, \mbox{ \rm  $\eta$-surrounds } G 
   \, \mbox { \rm \underbar{or} } \; B[0, \tau_{\exp}] \cap \corG = \emptyset
 \, \cond \, \cA_G   \Big ) > \clll
 \; \; \:
& \; & \;
\end{eqnarray}
where the conditioning is on the $\sigma$-field $\cA_G$
generated by $I, \,  B^{(1)}, \, B^{(3)} \one_{\{I=3 \} } \,$ and $\, \texp$. 
\end{pr}

\noindent{\sc Proof:}
  \, On the event $\{I=1\}$, the set $B[0, \tau_{\exp}]$
    is disjoint from $\corG$. \newline
To handle the case $\{I=2\}$, we use the strong Markov property at
time $\tauc$ and apply the estimate (\ref{eq:fast}) to $B^{(2)}$. 
Denoting 
$\taucJ= \min \{ t \geq \tau_c \, : \, B(t) \in \partial G \}$, this gives
\begin{eqnarray*}
& \P_0 \Big( B[\tauc,\taucJ] \;  
\mbox{ \rm $\eta$-surrounds } G  \: \mbox{ \rm and }
\taucJ < \tau_{\exp} < \taucG \cond I \geq 2 ; \, B^{(1)} \Big) & \\[3ex]
& \geq 
  \cfast(\eta) \P_0 \Big( \tauc < \tau_{\exp} < \taucG)
 \cond  I \geq 2 ; \, B^{(1)} \Big)   \, . &
\end{eqnarray*}
This proves (\ref{eq gosper-pf}) on the event $I=2$. 

Only the case $I=3$ remains.
 By using the strong Markov property at time $\tauc$ and
 applying Lemma \ref{lem 3 parts} to $B^{(2)}$, we see that
 for any $A \subset \partial G$,
$$
\P_0 \Big(B[\tau_C,\taucG] \; \,
   \mbox{ \rm $\eta$-surrounds } G \; \mbox{\rm and } 
  B(\taucG) \in A \cond B^{(1)} \Big)
 \geq \clool(\eta) \P_0 \Big(B(\taucG) \in A \cond B^{(1)} \Big) \, .
$$
In other words,
$$
\P_x \Big(  B[\tau_C,\taucG] \; \,
   \mbox{ \rm $\eta$-surrounds } G \cond B^{(1)}, \, B(\taucG) \Big)
 \geq \clool(\eta) \, .
$$ 
An application of the strong Markov property 
at time $\taucG$ shows that this lower bound
is still valid if we insert an additional conditioning on 
$I \geq 2$ and on $B^{(3)}$. \newline
Finally, since $\P_0(I =3 \| I \geq 2) \geq e^{-1}/2$ and 
$\texp$ is conditionally independent of $B[0,\taucG]$
given $I=3$, this completes the proof of the proposition.
$\Cox$

To obtain the uniformity in Theorem~\ref{th main}, we will need the
following general observation.
\begin{lem} \label{lem semicont}
Let $\Gamma : [0,\infty) \rightarrow \complex$ be any continuous path
and let  $t > 0$. For any open disk $U$ intersecting
$\fr(\Gamma[0,t])$ such that  
$\Gamma (t) \notin \Ubar $, 
there is a $\delta > 0$ such that for any $s \in
[t , t + \delta ]$, we have 
\be \label{eq:contain}
U \cap \fr ( \Gamma [0,t]) \, \supset \,
 U \cap \fr (\Gamma [0,s]) \, \neq \emptyset \, .
\ee
\end{lem}
\noindent{\sc Proof:} 
By hypothesis $U$ intersects the unbounded component, $\Omega$.
of $\Gamma [0,t]^c$, so there is a point $u\in U$ with and
an unbounded curve starting from $u$ and contained in $\Omega$.
Using the convexity of $U$, we can append to this curve a line-segment
connecting $u$ to a nearest point $x$ on $\Gamma [0,t]$,
and thus obtain an unbounded curve $\psi$
starting at $x$ and contained in $\Omega \cup \{x\}$.
Choose $\delta>0$ small enough so that 
$\Gamma [t,t + \delta]$ is disjoint from the curve $\psi$.
This gives the right-hand side of (\ref{eq:contain}) 
for $s \in [t , t + \delta ]$.
If we also require that $\Gamma [t,t + \delta]$ is disjoint from 
$U$,  then the left-hand side of (\ref{eq:contain})
follows from the general fact that
$$
\fr ( \Gamma [0,s]) \subset \fr (\Gamma [0,t]) \cup \Gamma [t,s] \, .
$$
$\Cox$

\noindent{\sc Proof of Theorem~\ref{th main}:} $\:$
The random set $B[0, \tau_{\exp}] \setminus G^\circ$ is completely
determined by the variables generating the $\sigma$-field
$\cA_G$ defined in Proposition \ref{pr:brown}, so the
proposition implies that $K=B[0, \tau_{\exp}]$ satisfies
the hypothesis (\ref{eq gosper-wiggly}) of Theorem \ref{th any K}.

Since $B[0,\tau_{\exp}]$ has the same distribution as
$\sqrt{\tau_{\exp}} \cdot B[0,1]$,
this establishes the first assertion of
Theorem~\ref{th main}. 

 For $t>0$, Let $A_t$ be the event that 
$\dim (\fr (B[0,s]) \cap V) \geq 1 + \epsilon$ simultaneously for
all open disks $V$ that intersect $\fr(B[0,t])$ and have 
rational centers and radii. Theorem \ref{th any K}
and Proposition \ref{pr:brown} give $\P_0(A_1)>0$, and we must 
 show that 
${ \displaystyle \P_0  ( \cap_{t>0} A_t ) = 1}$.
Denote by $A_\Q= \cap_{s \in \Q_+} A_s$ the 
intersection over all positive rational times.  
Now
Brownian scaling and countable additivity 
 imply that $\P_0 (A_\Q)=1$, so it suffices to prove that
$A_\Q \subset A_t$ for all $t>0$.
Fix $t>0$ and an open disk $V$ that intersects $\frKBt$.
Since $\fr (B[0,t])$ is connected, it must intersect some 
(random) open disk $U = U(V,t)$ with rational center
and radius such that $U \subset V$ and $B(t) \notin \Ubar$.  By the 
previous lemma, 
there is a rational $s$ such that
$$
\fr (B[0,t]) \cap U \, \supset \, \fr (B[0,s]) \cap U \neq \emptyset \, .
$$
This implies that $A_\Q \subset A_t$,
and completes the proof of the theorem.
$\Cox$

Finally, we consider the {\bf planar Brownian bridge} $\Bbridge$,
which may be defined either by conditioning the Brownian path
to return to the origin, or by  $\Bbridge(t) = B(t)-tB(1)$ for $t \in [0,1]$.
For every $t<1$, the restrictions $\Bbridge|_{[0,t]}$
and $B|_{[0,t]}$ have mutually absolutely continuous laws
(these laws are measures on  the space of 
continuous maps from $[0,t]$ to the plane.)
Therefore by Theorem \ref{th main},
for every {\em fixed} $t \in (0,1)$,
\be \label{eq:bridge1}
\dim \Big(\fr(\Bbridge[0,t])  \Big) \, \geq 1 + \eps \; \, \mbox{ a.s.}
\ee
Consider a sequence of annuli $\{A_n\}$ of modulus $2^{-n}$
around the origin. The probability that $\Bbridge$
surrounds  the origin in $A_n$ is bounded away from 0,
so the Blumenthal 0--1 law implies that 
with probability 1, there is some rational
$t<1$ such that $\fr(\Bbridge[0,1]) = \fr (\Bbridge[0,t]) $
(see Burdzy and Lawler (1990)). 
Thus by (\ref{eq:bridge1}), with probability 1,
$$
\dim \Big(\fr(\Bbridge[0,1])  \Big) \, \geq
 \inf_{t \in \Q \cap (0,1)} \, \dim \Big(\fr(\Bbridge[0,t])  \Big)
 \geq \, 1 + \eps \,.
$$
\section{Concluding remarks}
It can be shown (Krzysztof Burdzy, personal communication)
 that $\dim (\frKB)$ is almost surely constant; 
this fact is not required for the
arguments in this paper.
The conjecture that the Brownian frontier has dimension $4/3$ is
related to well-known conjectures
concerning self-avoiding random walks, which in turn are a model for
long polymer chains. In that context, the exponent $4/3$
first appeared in the non-rigorous
considerations of Flory (1949); see also de Gennes (1991).


 Theorem \ref{th any K} is stated for general random sets,
rather than just Brownian motion, in view of potential
applications to the ranges and level-sets of other
stochastic processes.
%
Besides the range of Brownian motion, another natural random set that satisfies
the hypothesis of Theorem \ref{th any K} is the support
 of {\bf super-Brownian motion}, i.e. the intersection of all closed sets
that are assigned full measure by this measure-valued diffusion
throughout its lifetime. (For the definitions see, e.g.,
  Dawson, Iscoe, and Perkins (1989).)
Equivalently, this random set may be characterized as
the set of points ever visited by the path-valued process 
constructed by Le-Gall (1993). (This process is often referred to as 
``The Brownian snake''.).
We are grateful to Steve Evans for enlightening discussions of 
super-Brownian motion.

To allow for further applications,
 we state below a variant of Theorem \ref{th any K}
which obtains the same conclusions under weaker
hypotheses on the random set $K$.
We omit the proof, which requires the  estimates
obtained  by Pemantle (1994) for the probability that a Wiener sausage 
covers a straight line segment.

For any set $S \subset \complex$ and any $\eps > 0$, 
let $S^\eps$ denote the set $\{ x : |x-y| \leq \eps \mbox{ for some }
y \in S \}$.
  Say that $K$ {\em is $\eta , \delta$-flat} inside $S$ 
if there is some line segment $\ell$ of length $\eta \diam (S)$ covered by 
$K^{\delta \diam (S)}$, having $\ell^{\eta \diam (S)}$ inside $G$
with $\ell$ not topologically surrounded by $K \cap G \cap 
(\ell^{\delta \diam (S)})^c$.  
%
\begin{th} \label{th any K2}
Let $G_0$ be the Gosper Island, and
let $K$ be a random compact connected
subset of the plane.  
Suppose that for some $\delta_0 > 0$, the following hypothesis 
on $K$ is satisfied, where the supremum is over $r \in (0,1)$
and $\xx$ in the plane.
\begin{equation} \label{eq sup condition}
   \sup_{G = \xx + r G_0} \esssup 
   \P \left [ K \mbox{ is } \eta_0 , \delta_0 \mbox{-flat inside } 
   G \| K \cap G \neq \emptyset , \sigma (K \cap G^c) \right ] < 1 .
\end{equation}
Then there is an $\eps > 0$ for which
$\; \dim (\frK) \geq 1 + \eps \, $
with probability 1. 
\end{th}

\noindent{\sc Remark:} The intuition behind the two-part
definition of $\eta , \delta$-flatness is 
that for 
$\frK$ to be close to straight (thus for $K$ to be flat),
$K$ itself must nearly cover a line segment and this must 
happen somewhere that is not completely encircled by $K$.

For the special case when $K$ is the range of planar Brownian motion,
 it seems likely that methods directly adapted to this case will yield
better  estimates for $\dim (\frK)$
than those obtainable by our methods. Indeed, Gregory Lawler
has informed us that immediately
after he learned of our Theorem \ref{th main}
 (but without seeing its proof),  he proved
(using completely different methods)   
that the dimension of the Brownian frontier 
can be expressed in terms of the
 ``double disconnection exponent'' of Brownian motion.
This allowed Lawler to deduce that $\dim (\frKB) > 1.01$ a.s.,
by invoking recent estimates of Werner (1994) on disconnection exponents.
We refer the reader to Lawler's forthcoming paper
for this and several other striking results on the Brownian frontier.

Finally, we note an application to simple random walk on the
square lattice $\ZZ^2$.
Given a  subset $S$ of  $\ZZ^2$,
say that a lattice point $x \in S$ is on the {\bf outer boundary}
of $S$ if $x$ is adjacent to some point in the unbounded component
of $\ZZ^2 \setminus S.$
We remark that using the strong approximation
results of Auer (1990) and our construction of the Whitney
tree in Section 3, it is easy to derive the following.
\begin{cor} \label{cor:srw}
Let $\{S(k)\}$ denote simple random walk on $\ZZ^2$,
and let $\epsilon>0$ be as in Theorem \ref{th main}.
Then for every $\epsilon_1  < \epsilon$ we have
$$
\lim_{n \to \infty} 
\P\{ \mbox{There are more than } \: n^{(1+\epsilon_1)/2} \,
     \mbox{ points on the outer boundary of } \: S[0,n]. \} = 1 
$$
\end{cor}



\renewcommand{\baselinestretch}{1.0}

\noindent Christopher J. Bishop \\ Department
of Mathematics, SUNY at Stony Brook, Stony Brook, NY 11794-3651.

\noindent Peter W. Jones \\  Department  of Mathematics, Hillhouse Ave.,
Yale University,
New Haven, CT 06520 .

\noindent Robin Pemantle \\ Department 
of Mathematics, University of Wisconsin-Madison, Van Vleck Hall, 480 Lincoln
Drive, Madison, WI 53706 . 

\noindent Yuval Peres \\ Department of
Statistics, 367 Evans Hall University of California, Berkeley, CA 94720-3860.

\end{document}